\documentclass[12pt]{amsart}
\usepackage{amssymb}
\usepackage{graphicx}
\numberwithin{equation}{section}
\newtheorem{thm}{Theorem}[section]
\newtheorem{lem}[thm]{Lemma}
\newtheorem{conj}[thm]{Conjecture}
\newtheorem{prop}[thm]{Proposition}
\newtheorem{cor}[thm]{Corollary}

\newtheorem{rem}[thm]{Remark}

\newcommand{\lb}{\left[}
\newcommand{\rb}{\right]}
\newcommand{\lp}{\left\{}
\newcommand{\rp}{\right\}}
\newcommand{\ov}{\overline}

\newcommand{\rh}{\ti{H}}

\def\sgn{\operatorname{sgn}}

\def\im{\operatorname{im}}

\def\RR{{\mathcal{R}}}

\def\Ss{\mathfrak S}
\def\C{{\mathbb{C}}}

\def\Z{{\mathbb{Z}}}

\def\Q{{\mathbb{Q}}}

\def\no={\,{\,|\!\!\!\!\!=\,\,}}

\def\({\left(}
\def\){\right)}
\def\wt{\widetilde}

\def\sgn{\operatorname{sgn}}

\def\rh{\tilde H}

\def\beq{\begin{eqnarray}}
\def\eeq{\end{eqnarray}}
\def\beqn{\begin{equation} \begin{split}}
\def\eeqn{end{equation} \end{split}}
\begin{document}

\title[Torsion]{Torsion in the Matching Complex and  Chessboard Complex}
\author[Shareshian]{John Shareshian$^1$}
\address{Department of Mathematics, Washington University, St. Louis, MO}
\footnotetext[1]{Supported in part by NSF Grants
DMS 0070757 and DMS 0300483}
\email{shareshi@math.wustl.edu}

\author[Wachs]{Michelle L. Wachs$^2$}
\address{Department of Mathematics, University of Miami, Coral Gables, FL 33124}
\email{wachs@math.miami.edu}
\footnotetext[2]{Supported in part by NSF Grants
DMS 0073760 and DMS 0302310}

\date{\today \\ \small MR Subject Classifications: 05E25,  05E10, 55U10}

\begin{abstract}
Topological properties of the  matching complex were first studied
by Bouc in connection with Quillen complexes, and  topological
properties of the chessboard complex were first studied by Garst
in connection with Tits coset complexes.  Bj\"orner, Lov\'asz,
Vr\'ecica and {\v Z}ivaljevi\'c established bounds on the
connectivity of these  complexes and conjectured that these bounds
are sharp.  In this paper we show that the conjecture is true by
establishing the nonvanishing of
 integral homology in the degrees given by these bounds.  Moreover,  we show that for sufficiently large  $n$,
the bottom nonvanishing  homology of the matching complex
$M_n$ is an elementary 3-group, improving a result of Bouc, and that the bottom nonvanishing
homology of the chessboard complex $M_{n,n}$ is a 3-group of exponent at most 9.  When $n
\equiv 2 \bmod 3$, the bottom nonvanishing homology of $M_{n,n}$ is shown to be $\Z_3$.
 Our proofs rely on computer calculations, long exact sequences, representation
theory, and tableau combinatorics.
\end{abstract}

\maketitle

\tableofcontents

\section{Introduction}

A {\it matching} is a graph in which each vertex  is
contained in at most one edge. Given a graph $G=(V,E)$,  the collection of all
subgraphs $(V,F)$ of $G$ that are matchings forms an abstract simplicial
complex $M(G)$.  The vertices of $M(G)$ are
the edges of $G$, and the $k$-dimensional faces of $M(G)$ are the edge sets
$F$ of size $k+1$ such that
$(V,F)$
 is a matching.  If $G$ is the complete graph on
vertex set $[n] :=
\{1,2,\dots,n\}$, then
we write $M_n$ for $M(G)$. Similarly, if
$G$ is the complete bipartite graph with parts $[m]$ and
$[n]^\prime :=
\{1^\prime,2^\prime,\dots,n^\prime\}$
then we write $M_{m,n}$ for $M(G)$.

The complex $M_n$ is  called the {\it matching complex} and the
complex $M_{m,n}$ is called the {\em chessboard complex}.   A
piece of $M_7$ (taken from \cite{Bo})  is given in Figure~1.1
below. Here and throughout the paper, the vertex  of $M(G)$
labelled $ij$ represents the edge $\{i,j\}$ of the graph $G$. Each
$k$-dimensional face of the chessboard complex $M_{m,n}$
corresponds to a placement of $k+1$ nontaking rooks on an $m
\times n$ chessboard.  Indeed, a rook in the $i$th row and $j$th
column  corresponds to the edge $\{i,j^\prime\}$ in the bipartite
graph, which corresponds to the vertex $ij^\prime$ in $M_{m,n}$.
It is for this reason that the name ``chessboard complex'' is
used.

\vspace{-.6in}
\begin{center}
\includegraphics[width=4.7in]{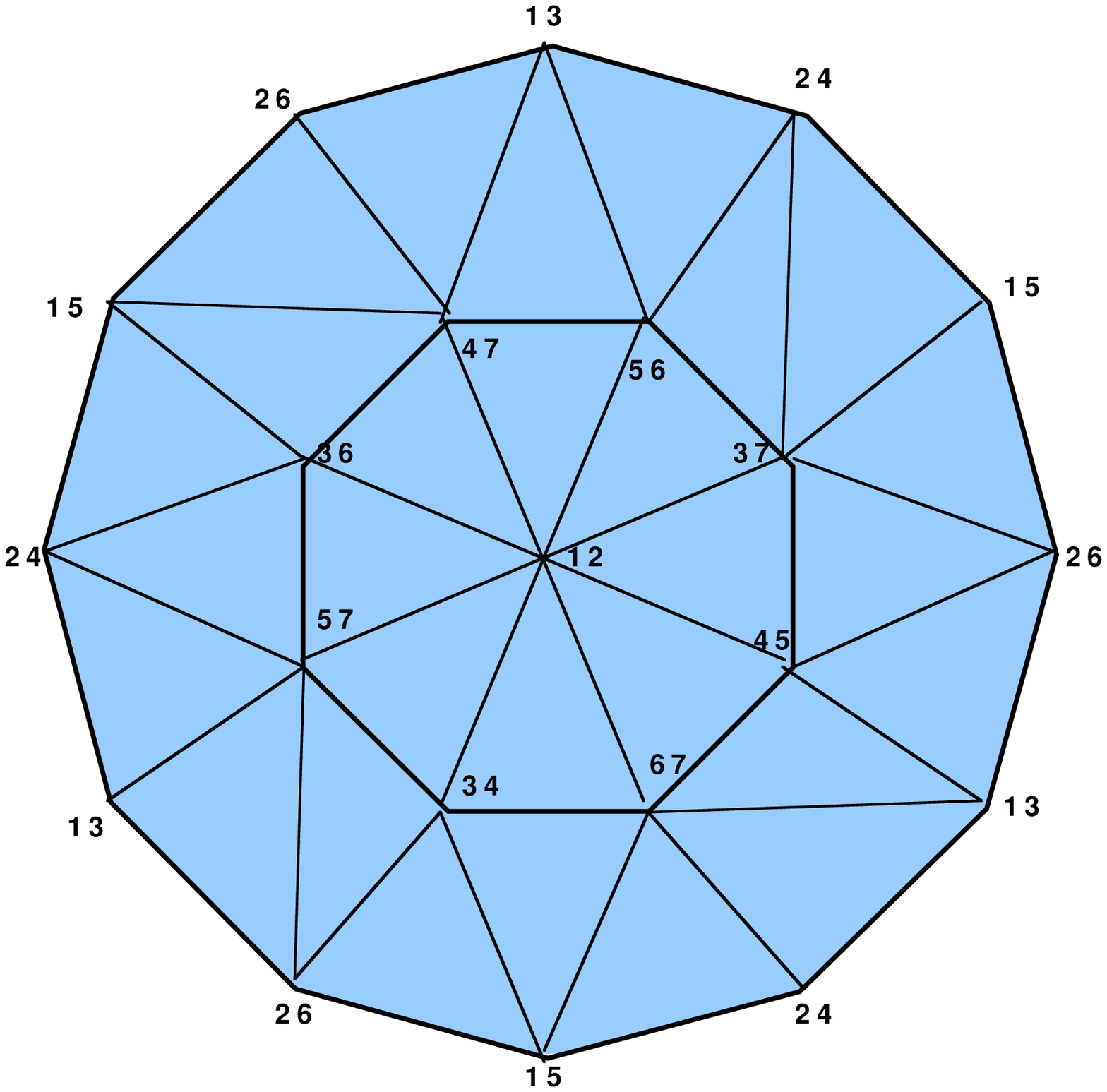}
\end{center}

\vspace{-.1in}\begin{center}{\bf Figure 1.1:} Piece of matching complex $M_7$ \end{center}

\vspace{.2in}

The matching complex, the chessboard complex and variations  have arisen in a variety of
fields such as group theory, representation theory,  commutative algebra, Lie theory,
computational geometry, and combinatorics;  see the survey article \cite{W1} and its
references.   Topological properties
of the matching complex were first studied by  Bouc \cite{Bo}, in connection
with the Quillen complex at the prime $2$ for the symmetric group.  Bouc
obtains several beautiful results.  He considers the representation of the symmetric group
$\mathfrak S_n$ acting on the  homology (over $\C$) of the
matching complex
$M_n$ and obtains a  decomposition into irreducibles. This yields a formula for the Betti
numbers in terms of standard Young tableaux.  Bouc also obtains  results on torsion in integral
homology, which we improve and extend to the chessboard complex in this paper.

Prior to Bouc's study of the matching complex,  the chessboard complex was
introduced in the 1979 thesis of Garst \cite{G}
dealing with Tits coset complexes.  Garst shows that for $m \le n$, $M_{m,n}$ is
Cohen-Macaulay  if and only if $2m-1 \le n$.  Garst also obtains  a
decomposition of the representation of
$\mathfrak S_n$ acting on the {\em top} homology of $M_{m,n}$ into
irreducibles, for  $m \le n$.   This computation is a precursor
of  Friedman and Hanlon's
\cite{FH} decomposition of the representation of $\Ss_m\times \Ss_n$ on {\em each}  homology
of $M_{m,n}$ into irreducibles.

 Questions on connectivity of the chessboard complex
were raised by \u  Zivaljevi\'c   and Vr\'ecica \cite{ZV} in
connection with some problems in computational geometry.   In
response to these questions, Bj\"orner, Lov\'asz, Vr\'ecica, {\v
Z}ivaljevi\'c \cite{BLVZ} obtained  bounds on connectivity of the
chessboard complex and the matching complex which are given in the
following theorem.  The bound for the matching complex is also an
immediate consequence of results in Bouc \cite{Bo}.

\begin{thm}[Bj\"orner, Lov\'asz,
Vr\'ecica, {\v Z}ivaljevi\'c \cite{BLVZ}, Bouc \cite{Bo}]
\label{conn} For positive integers $m,n$, let $$\nu_n = \lfloor
\frac{n+1} 3 \rfloor -1 \quad \text{ and } \quad \nu_{m,n} =
\min\{m,n, \lfloor \frac {m+n+1} 3 \rfloor\} -1.$$ Then the
matching complex $M_n$ is $(\nu_n-1)$-connected and the chessboard
complex $M_{m,n}$ is $(\nu_{m,n}-1)$-connected. Consequently, for
all $t < \nu_n$, \beq\label{blvzeq1}\rh_t(M_n) = 0, \eeq and for
all $t < \nu_{m,n}$, \beq \label{blvzeq2}\rh_t(M_{m,n}) = 0. \eeq

\end{thm}

\begin{rem}{\rm Throughout this paper, by homology of a simplicial complex $\Delta$, we mean
reduced simplicial homology $\rh_*(\Delta)$ over the integers,  unless otherwise stated.}
\end{rem}

 It is
conjectured in
\cite{BLVZ} that the connectivity bounds of  Theorem~\ref{conn} are sharp.  The $n \equiv 0,1
\bmod 3$ cases  of the conjecture for the matching complex had already been  established by
Bouc \cite{Bo} who  proved the following result.

\begin{thm}[Bouc\cite{Bo}] \quad \label{Btor}
\begin{itemize}
\item[(i)] $\rh_{\nu_{n}}(M_{n})$ is finite if and only if   $n
\ge 7$ and $ n \notin \{8,9,11\}$.
\item[(ii)] If $n \equiv 1 \bmod 3$ and $n \ge 7$ then $\tilde H_{v_n}(M_n) \cong
\Z_3$.
\item[(iii)] If $n \equiv 0 \bmod 3$ and $n \ge 12$ then $\tilde H_{v_n}(M_n)$ is a
nontrivial $3$-group  of exponent at most
$9$.
\end{itemize}
\end{thm}

\begin{rem}{\rm Statement (i) is not explicitly stated in \cite{Bo}, but follows easily from
the formula for the Betti numbers given in \cite{Bo}.}
\end{rem}

 One can  see the   $3$-torsion in $\rh_1(M_7)$ by
looking at Figure~1.1.  The union  of the triangles shown  is bounded by  $3 z$ where $$z
= (13,24)+(24,15)+(15,26)+(26,13).$$  Bouc shows that $ z$
is not a boundary; so $z$ is a $3$-torsion element.

 Friedman and Hanlon
\cite{FH}  derive  the following analogy of Theorem~\ref{Btor}~(i), which settles the
chessboard complex version of the conjecture in the case that $n > 2m-5$, but
leaves the conjecture unresolved  the case that $m \le n \le 2m-5$.  Their result is  a
consequence of their formula  for  the Betti numbers of the chessboard
complex  derived in
\cite{FH} (see Theorem~\ref{FH2}).

\begin{thm}[Friedman and Hanlon \cite{FH}] \label{FH} Let $m \le n$. Then the group
$\rh_{\nu_{m,n}}(M_{m,n})$ is finite if and only if  $ n
\le 2m-5$ and
$(m,n) \notin \{(6,6),(7,7),(8,9)\}$.
\end{thm}

In this paper we pick up where Bouc and Friedman-Hanlon left off.
We prove the Bj\"orner-Lov\'asz-Vr\'ecica-{\v Z}ivaljevi\'c
conjecture   in the  cases that were left unresolved in Bouc's
work and Friedman-Hanlon's work (see Theorem~\ref{conj}).
Moreover, we  prove the following result which improves Theorem~\ref{Btor}  by handling  the remaining $n \equiv 2 \bmod 3$ case
and making the exponent precise in all cases.

\begin{thm} \label{mtor} For  $n = 7,10$ or $n \ge 12$ (except possibly $n=14 $),
$ \tilde H_{\nu_n}(M_n)$ is a nontrivial elementary $3$-group.
\end{thm}

We also prove the following analogous result for the chessboard complex.

\begin{thm}  Let $m \le n$. \label{ctor}
\begin{itemize}
\item[(i)] If $m+n \equiv 1 \bmod 3$ and $ n \le 2m-5$ then $\tilde H_{v_{m,n}}(M_{m,n}) \cong
\Z_3$.
\item[(ii)] If $m+n \equiv 0 \bmod 3$ and $ n \le 2m-9$ then
$\tilde H_{v_{m,n}}(M_{m,n})$ is a  nontrivial $3$-group  of
exponent at most
$9$.
\item[(iii)] If $m+n \equiv 2 \bmod 3$ and $ n \le 2m-13$ then $\tilde H_{v_{m,n}}(M_{m,n})$
is a  nontrivial $3$-group  of exponent at most
$9$.
\end{itemize}
\end{thm}

Bouc proves the $1 \bmod 3$ case of Theorem~\ref{Btor}
 using induction.  His main tool is a  long exact sequence which provides the induction step
and also enables him to derive  the $0 \bmod 3$ case from  the $1 \bmod 3$ case. Bouc's
``hand'' calculation of
$\rh_{\nu_7}(M_7)$ provides the base step of the induction.    Here we further exploit Bouc's
long exact sequence to derive the $2 \bmod 3$ case from the $0 \bmod 3$ case, and we use a
computer calculation to provide another base case
$\rh_{\nu_{12}}(M_{12})$ which enables us to bring the exponent down to $3$ in
Theorem~\ref{mtor}.

The proof of Theorem~\ref{ctor}, while patterned on the proof of the Theorem~\ref{mtor}, is
much more difficult.  An essential ingredient is  an interesting basis for the
{\it top} homology of the chessboard complex.  The construction of this basis has a
surprising reliance  on a result in tableau combinatorics, namely the classical
Robinson-Schensted correspondence.

The computer program that we use for computing homology in the base steps,  was first
developed by  Heckenbach and later improved by Dumas, Heckenbach, Saunders and Welker
\cite{DHSW}.  With this software, one can produce the following tables.
$$\begin{array}{|c||c|}
\hline  n  \,  &\rh_{\nu_n}(M_n) \\ \hline \hline
  2 &  0  \\ \hline
 3 &  \Z^2  \\ \hline
 4 &   \Z^2  \\ \hline
 5 &    \Z^6  \\ \hline
 6 &   \Z^{16} \\ \hline
 7 &    \Z_3 \\ \hline
 8 &   \Z^{132}\\ \hline
 9 &   \Z^{42}\oplus \Z_3^8   \\ \hline
10 &    \Z_3  \\ \hline
11 &   \Z^{1188}\oplus Z_3^{45} \\ \hline
12 &  \Z_3^{56}
\\ \hline
13 &  \Z_3  \\ \hline
14 & ? \\ \hline
\end{array}$$
\begin{center} {\bf Table 1.1:} Bottom nonvanishing homology
$\rh_{\nu_n}(M_{n})$
\end{center}

\vspace{.1in}

$$ \begin{array}{|c||c|c|c|c|c|c|c|}
\hline m  \backslash n \,  & 2  & 3  & 4 & 5 &
6 & 7 & 8\\ \hline\hline
 2 &  \Z & \Z & \Z^5 & \Z^{11} & \Z^{19} & \Z^{29} & \Z^{41} \\ \hline
 3 &   & \Z^4 & \Z^2 &\Z^{14} &\Z^{47} & \Z^{104} & \Z^{191} \\ \hline
 4 &   &  & \Z^{15} & \Z^{20}\, & \Z^5 & \Z^{225} & \Z^{641} \\ \hline
 5 &   &  &  &  \Z_3 & \Z^{152} & \Z^{98} &\Z^{14} \\ \hline
 6 &   &  &  &  & \Z^{25}\oplus  \Z_3^{10} \,\,\,&\Z_3 & \Z^{1316}\\ \hline
 7 &   &  &  &  &  &  \Z^{588} \oplus \Z_3^{66} & ?\\ \hline
\end{array}$$
\begin{center} {\bf Table 1.2:} Bottom nonvanishing homology  $\rh_{\nu_{m,n}}(M_{m,n})$
\end{center}

\vspace{.1in}

Unfortunately we have not been able to get output for $n \ge 14$ nor for $m \ge 7$ and $n \ge
8$.
 This is what is responsible for the  gap at $n = 14$ in Theorem~\ref{mtor} and the
lack of precision with respect to the exponent in Theorem~\ref{ctor}.  Indeed,  in
Theorems~\ref{conj99} and \ref{genmain}, we show that if we could determine   the exponent
of the Sylow 3-subgroup of
$\tilde H_{v_{7,8}}(M_{7,8})$   or  the exponent of $\tilde H_{v_{9,9}}(M_{9,9})$ to be 3,
then we could conclude that the exponent of
$\tilde H_{v_{m,n}}(M_{m,n})$ is 3 for all
$m,n$ that satisfy the conditions of Theorem~\ref{ctor}.

The paper is organized as follows.  In Section~\ref{les}, notation
is established and the long exact sequences are derived. In
Section~\ref{BLVZsec}, we prove  the
Bj\"orner-Lov\'asz-Vr\'ecica-{\v Z}ivaljevi\'c  connectivity
conjecture.  The torsion result for the matching complex,
Theorem~\ref{mtor}, is proved in  Section~\ref{mtorsec}.

Sections~\ref{ctorsec}, \ref{top}, \ref{infinsec} and \ref{subsec} are devoted to the
chessboard complex.   The proof of Theorem~\ref{ctor} is given in Section~\ref{ctorsec}.
Partial results on torsion in the finite groups
$\tilde H_{v_{m,n}}(M_{m,n})$ not covered by Theorem~\ref{ctor} can also be found in
Section~\ref{ctorsec}.  The basis for the top homology of the chessboard complex used in theproof of Theorem~\ref{ctor} is constructed in Section~\ref{top}.

In Section  \ref{infinsec}, we deal with torsion in the case of infinite  $\tilde
H_{v_{m,n}}(M_{m,n})$.  Here we use Friedman and Hanlon's representation theoretic
result to show that $\rh_{\nu{m,n}}(M_{m,n})$ is torsion-free when $n = 2m-2$.  We give
 conjectures for the other cases of infinite homology.

In Section~\ref{subsec}, we discuss the subcomplex of the square
chessboard complex
$M_{n,n}$ obtained by deleting a diagonal from the chessboard.  This complex was shown to be
$(\nu_{2n}-1)$-connected by Bj\"orner
and Welker
\cite{BW} as a consequence of a more general result of Ziegler \cite{Z} on nonrectangular
boards.  Here we
show that the Bj\"orner-Welker-Ziegler bound is sharp.

In Section~\ref{shellsec}, we  answer another question of
Bj\"orner, Lov\'asz, Vr\'ecica, and {\v Z}ivaljevi\'c \cite{BLVZ}.
Given the connectivity bounds on $M_n$ and $M_{m,n}$, they ask
whether the $\nu_n$-skeleton of $M_n$ and the $\nu_{m,n}$-skeleton
of $M_{m,n}$ are shellable.       Ziegler \cite{Z} answers this
question affirmatively for the chessboard complex by establishing
vertex decomposability.  In Section 10, we answer the question
affirmatively for  the matching complex. We remark that in
subsequent work, Athanasiadis \cite{A} improves this result by
establishing vertex decomposability.

In Section~\ref{ranksec}, bounds on the ranks of the finite 3-groups $ \rh_{\nu_n}(M_n)$ and
$\tilde H_{\nu_{m,n}}(M_{m,n})$ are derived.  This extends bounds given by Bouc for the $n
\equiv 0,1 \bmod 3$ cases of the matching complex.

\section{Bouc's long exact sequence} \label{les}

In \cite{Bo}, Bouc produces a long exact sequence which enables
him to prove that $\rh_t(M_n)=0$ for $t<\nu(n)$ and to obtain
Theorem~\ref{Btor}.  As we will see in Section~\ref{BLVZsec}, it
is easy to use Bouc's sequence to show that $\rh_{\nu_n}(M_n) \neq
0$ when $n \equiv 2 \bmod 3$, thereby establishing the matching
complex case of the Bj\"orner-Lov\'asz-Vr\'ecica-{\v Z}ivaljevi\'c
conjecture. This sequence  will also play a role in the proof of
Theorem~\ref{mtor} given in Section~\ref{mtorsec}.   In this
section, we  present Bouc's long exact sequence and an analogous
sequence for the chessboard complex. The analogous sequence will
be used to prove the chessboard complex version of the
Bj\"orner-Lov\'asz-Vr\'ecica-{\v Z}ivaljevi\'c conjecture in
Section~\ref{BLVZsec},  and to prove Theorem \ref{ctor} in
Section~\ref{ctorsec}.

We use  standard notation,  $(C_*(\Delta),\partial)$ and $Z_*(\Delta)$, for the chain complex
and  the cycle group, respectively,  of a simplicial complex $\Delta$.  For $z \in Z_*(\Delta)$,
we let
$\bar z$ denote the homology class of $z$ in $\tilde H_*(\Delta)$.

\subsection{The long exact sequence for $M_n$} \label{les1}

In order to state Bouc's result in a manner that will be useful to us, we must introduce
some additional notation.  For finite set $A$, let $M_A$ be the matching complex on the
complete graph with vertex set $A$.

For disjoint subsets $A,B \subseteq [n]$, if $z_1$ and $z_2$ are oriented simplices of
$M_A$ and $M_B$, respectively, then $z_1 \wedge z_2$ will denote the oriented simplex of
$M_{A\cup B}$ obtained by concatenating $z_1$ and $z_2$.  We define a homomorphism
$$\bigwedge: C_{s-1}(M_A) \otimes C_{t-1}(M_B) \to C_{s+t-1}(M_{A
\cup  B})$$ by letting $z_1 \otimes z_2 \mapsto z_1 \wedge z_2$ for all oriented simplices
$z_1,z_2$.   This induces  a homomorphism
$$\bigwedge: \rh_{s-1}(M_A)
\otimes \rh_{t-1}(M_B)
\to \rh_{s+t-1}(M_{A
\cup  B}),$$ defined by $\ov {z_1} \wedge \ov {z_2} = \ov{z_1 \wedge z_2}$ for all
$z_1\in Z_{s-1}(M_A)$ and $z_2\in Z_{t-1}(M_B)$. (We  write
$z_1
\wedge z_2$ instead of $\bigwedge(z_1 \otimes z_2)$  and $\ov{z_1}
\wedge \ov{z_2}$ instead of $\bigwedge(\ov{z_1} \otimes \ov{z_2})$ and note that
$z_1 \wedge z_2$ is a cycle.)

For $a =1,2$ and $i = 3,\dots,n$, let
$$\phi_{a,i}: \rh_{t-1}(M_{[n]\setminus\{1,2,i\}}) \to \rh_t(M_n)$$ be the homomorphism
defined by
$$\phi_{a,i}(\ov{z}) =\ov{ai-12 } \,\, \wedge \,\, \ov{z} .$$ This determines the
homomorphism
$$\phi : \bigoplus_{\scriptsize\begin{array}{c} a\in \{1,2\}

\\ i\in [n]\setminus\{1,2\} \end{array} } \!\!\rh_{t-1}(M_{[n]\setminus\{1,2,i\}}) \, \, \to
\,\,
\rh_t(M_n),$$ defined by letting $\phi(\ov{z}) = \phi_{a,i}(\ov{z})$ for each
$\ov{z}$ in each
$(a,i)$-summand.

For $i \ne j \in \{ 3,\dots,n\}$, let
$$\psi_{i,j} : C_t(M_n) \to C_{t-2}(M_{[n]\setminus\{1,2,i,j\}})$$ be the map defined by
letting $$\psi_{i,j}(x) = \lp
\begin{array}{ll}  y   & {\rm if } \,\,\, x = 1i \wedge 2j \wedge y \mbox{ for some } y \in
C_{t-2}(M_{[n]\setminus\{1,2,i,j\}})
\\ 0 & {\rm otherwise,}
\end{array}
\right. $$ for each oriented simplex $x$. It is straightforward to show that the induced map
$$\psi_{i,j}:\rh_t(M_n) \to \rh_{t-2}(M_{[n]\setminus\{1,2,i,j\}})$$ given by
$\psi_{i,j}(\ov{z}) = \ov{\psi_{i,j}(z)}$ is a well-defined homomorphism as is the map
$$\psi:\rh_t(M_n) \to \bigoplus_{i\ne j \in
[n]\setminus\{1,2\}}\rh_{t-2}(M_{[n]\setminus\{1,2,i,j\}})$$ given by
$\psi(\ov{z}) = ({\psi_{i,j}(\ov z)}).$

For $a =1,2$,\quad  $h,i,j = 3,\dots,n$ and $i \ne j$, define
$$\delta^{i,j}_{a,h}:\rh_t(M_{[n]\setminus\{1,2,i,j\}}) \to
\rh_t(M_{[n]\setminus\{1,2,h\}})$$ by
$$\delta^{i,j}_{a,h}(\ov{z}) = \lp
\begin{array}{rl}
\ov{z} & {\rm if } \;\;a=1 \;\;{\rm and } \;\;h=i\\ -\ov{z} & {\rm if } \;\;a=2
\;\;{\rm and } \;\;h=j\\ 0 & {\rm otherwise,}
\end{array}
\right. $$ for $z \in Z_t(M_{[n]\setminus\{1,2,i,j\}})$.  Again it is straightforward to show
that
$\delta^{i,j}_{a,h}$ is a well-defined homomorphism as is  the homomorphism
$$\delta:\bigoplus_{i\ne j \in
[n]\setminus\{1,2\}} \!\!\rh_t(M_{[n]\setminus\{1,2,i,j\}}) \,\,\,\to
\bigoplus_{\scriptsize\begin{array}{c} a\in\{1,2\}
\\ h\in [n]\setminus\{1,2\} \end{array} } \!\!\rh_t(M_{[n]\setminus\{1,2,h\}})$$ defined by
letting
$\delta(\ov{z}) = (\delta^{i,j}_{a,h}(\ov{z}))$ for each
$\ov{z}$ in each
$(i,j)$-summand.

We can now state Bouc's result.  For the sake of completeness,
 we will include a proof.

\begin{lem}[{\cite[Lemma 9]{Bo}}] The sequence
\[
\cdots \stackrel{\delta}{\rightarrow}\bigoplus_{\scriptsize\begin{array}{c} a\in\{1,2\}
\\ h\in [n]\setminus\{1,2\} \end{array} } \!\!
\rh_{t-1}(M_{[n]\setminus\{1,2,h\}}) \stackrel{\phi}{\rightarrow}
\rh_t(M_n) \stackrel{\psi} {\rightarrow}
\]
\[\bigoplus_{i \ne j \in [n]\setminus\{1,2\}}
\rh_{t-2}(M_{[n]\setminus\{1,2,i,j\}})
\,\,\stackrel{\delta}{\rightarrow} \!\!
\bigoplus_{\scriptsize\begin{array}{c} a\in\{1,2\}
\\ h\in [n]\setminus\{1,2\} \end{array} } \!\!\rh_{t-2}(M_{[n]\setminus\{1,2,h\}})
\stackrel{\phi}{\rightarrow}\cdots\qquad\qquad\qquad\qquad\qquad\qquad\qquad
\]  is exact. \label{lesbouc}
\end{lem}

\begin{proof} For any graph $G$ on vertex set $[n]$, let $E(G)$ denote the edge set of $G$,
and for$v
\in [n]$, let
$N_G(v)$ denote the  set of neighbors of $v$, that is,
$$N_G(v) = \{u \in V : \{u,v\} \in E(G)\}.$$
Define
\[ X_n:=\lp G \in M_n : |(N_G(1) \cup N_G(2)) \setminus \lp 1,2 \rp|
\leq 1 \rp.
\] Then $X_n$ is a subcomplex of $M_n$, and we examine the standard long exact sequence
\[
\cdots \stackrel{\partial_\ast}{\rightarrow}\rh_t(X_n)
\stackrel{i_\ast}{\rightarrow} \rh_t(M_n)
\stackrel{\pi_\ast}{\rightarrow} \rh_t(M_n,X_n)
\stackrel{\partial_\ast}{\rightarrow}\rh_{t-1}(X_n)
\stackrel{i_\ast}{\rightarrow} \cdots
\] (see \cite[Theorem 23.3]{Mu}).

Let $P_n$ be the subcomplex of $X_n$ consisting  those $G \in X_n$ such that either
$\lp 1,2 \rp \in E(G)$ or both $1$ and $2$ are isolated in $G$. Since $P_n $ is a cone over
$M_{[n]\setminus \{1,2\}}$, it is  acyclic.  Hence the  natural projection of chain
complexes induces an isomorphism
\[
\tau:\rh_t(X_{n}) \rightarrow \rh_t(X_{n},P_{n}).
\]

For $a
\in \lp 1,2
\rp$ and
$h
\in
\lb n
\rb
\setminus
\lp 1,2
\rp$, let
$$\alpha_{a,h}: C_t(X_n,P_n) \to C_{t-1}(M_{[n]\setminus \{1,2,h\}})$$
be the map defined by letting
$$ \alpha_{a,h}(x) = \begin{cases} y &\mbox{if } x = ah \wedge y \mbox{ for some } y \in
C_{t-2}(M_{[n]\setminus\{1,2,h\}})\\ 0 &\mbox{otherwise,}
\end{cases}$$
for each oriented simplex $x$.
 It is straightforward to show that the induced map
$$\alpha_{a,h}:\rh_t(X_n,P_n) \to \rh_{t-1}(M_{[n]\setminus\{1,2,h\}})$$ given by
$\alpha_{i,j}(\ov{z}) = \ov{\alpha_{i,j}(z)}$, is a well-defined homomorphism as is the map
$$\alpha:\rh_t(X_n,P_n) \,\,\to \!
\bigoplus_{\scriptsize\begin{array}{c} a \in \{1,2\} \\ h
\in [n]\setminus\{1,2\}\end{array}}\!\rh_{t-1}(M_{[n]\setminus\{1,2,h\}})$$ given by
$\alpha(\ov{z}) = ({\alpha_{a,h}(\ov z)}).$ If we define
\[
\gamma_{a,h}:\rh_{t-1}(M_{\lb n \rb \setminus \lp 1,2,h \rp})
\rightarrow \rh_t(X_n,P_n)
\] by
\[
\ov{w} \mapsto \ov{a h \wedge w}
\] then
\[
\gamma:=\bigoplus_{a,h}\gamma_{a,h}
\] is a well-defined inverse for $\alpha$.  We now have an isomorphism
\[
\alpha\tau:\rh_t(X_n) \rightarrow
\bigoplus_{\scriptsize\begin{array}{c} a\in
\{1,2\}\\h \in [n]\setminus\{1,2\} \end{array}}\!\!\rh_{t-1}(M_{\lb n \rb \setminus \lp 1,2,h
\rp}).
\]

It is straightforward to show that  the map  $$ \beta_{i,j}:\rh_t(M_n,X_n) \rightarrow
\rh_{t-2}(M_{\lb n
\rb \setminus \lp 1,2,i,j \rp})$$
 induced by the restriction of $\psi_{i,j}$ to $C_t(M_n,X_n)$ is  a
well-defined homomorphism for all $i,j
\in
\lb n
\rb
\setminus
\lp 1,2
\rp$ with $i \ne j$.
 Define
\[
\beta:\rh_t(M_n,X_n) \rightarrow \bigoplus_{\scriptsize\begin{array}{c} i,j\in
[n]\setminus\{1,2\} \\i\ne j\end{array}}\!\!\rh_{t-2}(M_{\lb n
\rb \setminus \lp 1,2,i,j \rp})
\] by
\[
\ov{z} \mapsto ( \beta_{i,j}(\ov{z})).
\] If we define
\[
\mu_{i,j}:\rh_{t-2}(M_{\lb n \rb \setminus \lp 1,2,i,j \rp})
\rightarrow \rh_t(M_n,X_n)
\] by
\[
\ov{w} \mapsto \ov{1i \wedge 2j \wedge w}
\] then
\[
\mu:=\bigoplus_{i,j}\mu_{i,j}
\] is an inverse for $\beta$.  The result now follows from the fact that the diagram
\begin{equation}  \begin{split} \nonumber
\cdots \stackrel{\partial_\ast}{\longrightarrow} \qquad  \rh_t & (X_n)
\quad \stackrel{i_\ast}{\longrightarrow} \quad \rh_t(M_n) \quad
\stackrel{\pi_\ast}{\longrightarrow} \quad \rh_t(M_n,X_n) \quad
\stackrel{\partial_\ast}{\longrightarrow} \cdots
\\ \nonumber &
\big{\downarrow} \scriptsize{\alpha\tau} \qquad\qquad\qquad\quad
\big{\downarrow} \scriptsize{\mbox{id}} \qquad\qquad\qquad\qquad
\quad \big{\downarrow} \scriptsize{\beta}
\\ \nonumber
\cdots \stackrel{\delta}{\rightarrow} \bigoplus_{a,h}
\rh_{t-1}&(M_{[n]\setminus\{1,2,h\}}) \stackrel{\phi}{\rightarrow}
\rh_t(M_n) \stackrel{\psi}{\rightarrow} \bigoplus_{i,j}
\rh_{t-2}(M_{[n]\setminus\{1,2,i,j\}})
\stackrel{\delta}{\rightarrow} \cdots \end{split}\end{equation}
commutes.
\end{proof}

\subsection{The long exact sequence for $M_{m,n}$} \label{les2}

For any subset $$Y = \{y_1,y_2,\dots,y_k\} \subseteq [n],$$ let $$Y^\prime
:=\{y_1^\prime,y_2^\prime,\dots,y_k^\prime\} \subseteq [n]^\prime.$$   For
$X\subseteq [m]$ and $Y
\subseteq [n]$, let
$M_{X,Y}$ be the chessboard complex on
$X$ and
$Y^\prime$.  In other words, $M_{X,Y}$ is the matching complex on the complete bipartite graph
whose parts are $X$ and $Y^\prime$. Then
$M_{X,Y}$ is a subcomplex of the matching complex $M_{X \uplus Y^\prime}$, and the chain complex
$C_\ast(M_{X,Y})$ is embedded in the complex
$C_\ast(M_{X \uplus Y^\prime})$.

After appropriate changes in notation, restrictions of the
various functions defined in Section
\ref{les1} will be used to produce a long exact sequence similar to the one described in
Lemma \ref{lesbouc}.  In particular, if
$X=X_1 \biguplus X_2$ and $Y=Y_1 \biguplus Y_2$ then the restriction of the homomorphism
$\bigwedge$ defined in Section
\ref{les1} gives a homomorphism
\[
\bigwedge:\rh_{s-1}(M_{X_1,Y_1}) \otimes \rh_{t-1}(M_{X_2,Y_2})
\rightarrow \rh_{s+t-1}(M_{X,Y}).
\]
In Section \ref{les1}, the graph vertices $1,2$ were distinguished in order to produce the
desired long exact sequence.  Here, we distinguish the graph vertices $1,1^\prime$.
For  $ i   \in \lb m \rb\setminus\{1\}$, define
\[
\phi_{ i  }:\rh_{t-1}(M_{\lb m \rb \setminus \lp 1 , i   \rp,\lb n
\rb \setminus \lp  1  \rp}) \rightarrow \rh_t(M_{m,n})
\] by
\[
\ov{z} \mapsto (\ov{ 11^\prime   -  i1^\prime }) \wedge \ov{z},
\] and for
$ j  \in \lb n \rb \setminus \lp  1  \rp$, define
\[
\phi^\prime_{ j }:\rh_{t-1}(M_{\lb m \rb \setminus \lp 1 \rp,\lb n \rb
\setminus \lp  1 , j  \rp}) \rightarrow \rh_t(M_{m,n})
\] by
\[
\ov{z} \mapsto (\ov{1  1^\prime -1  j^\prime }) \wedge \ov{z}.
\] For ease of notation, we define
\begin{eqnarray*}
\rh_t(1) & := & \bigoplus_{ j \in \lb n \rb \setminus \lp  1
\rp}\rh_t(M_{\lb m \rb \setminus \lp 1 \rp,\lb n \rb \setminus
\lp  1 , j  \rp}), \\ \rh_t( 1^\prime ) & := & \bigoplus_{ i   \in \lb m
\rb \setminus \lp 1 \rp}\rh_t(M_{\lb m \rb \setminus \lp 1, i
\rp,\lb n \rb \setminus \lp  1  \rp}).
\end{eqnarray*} The maps $\phi_{ i }$ and $\phi^\prime_{ j}$ together determine a unique
homomorphism
\[
\phi:\rh_{t-1}(1^\prime) \oplus \rh_{t-1}( 1 ) \rightarrow
\rh_t(M_{m,n}).
\]

For $ i  \in \lb m \rb \setminus \lp  1  \rp$ and $ j   \in \lb n
\rb \setminus \lp 1 \rp$, define
\[
\psi_{i,j  }:C_t(M_{m,n}) \rightarrow C_{t-2}(M_{\lb m \rb
\setminus \lp 1, i   \rp,\lb n \rb \setminus \lp  1 , j  \rp})
\] by
\[ x \mapsto \lp
\begin{array}{ll}  y   & \mbox{if } x=1 j^\prime  \wedge  1^\prime  i \wedge  y \mbox{ for some } y \in
C_{t-2}(M_{\lb m \rb
\setminus \lp 1, i   \rp,\lb n \rb \setminus \lp  1 , j  \rp})
\\ 0 &
\mbox{otherwise.}
\end{array}
\right.
\] As in Section \ref{les1}, $\psi_{ i,j  }$ induces a homomorphism, also called
$\psi_{i,j}$, from $\rh_t(M_{m,n})$ to
$\rh_{t-2}(M_{\lb m \rb \setminus \lp 1, i   \rp,\lb n \rb
\setminus \lp  1 , j  \rp})$. We define
\[
\psi:\rh_{t}(M_{m,n}) \rightarrow \bigoplus_{\scriptsize{
\begin{array}{c}  i  \in \lb m \rb \setminus \lp  1  \rp \\  j   \in \lb
n
\rb \setminus \lp 1 \rp \end{array}}}\rh_{t-2}(M_{\lb m \rb
\setminus \lp 1, i   \rp,\lb n \rb \setminus \lp  1 , j  \rp}),
\]  by
\[
\ov{z} \mapsto (\psi_{ i,j  }(\ov{z})).
\]

For $ i  \in \lb m \rb \setminus \lp  1  \rp$ and $ j   \in \lb n
\rb \setminus \lp 1 \rp$ define
\begin{eqnarray*} \delta^{i, j }: \rh_t(M_{\lb m \rb  \setminus \lp 1,
i\rp,\lb n \rb
\setminus
\lp  1 , j  \rp})  \rightarrow & &\\  \rh_t(M_{\lb m
\rb \setminus \lp 1,i \rp,\lb n \rb \setminus \lp  1 \rp})
&\oplus &
\rh_t(M_{\lb m \rb \setminus \lp 1   \rp,\lb n
\rb \setminus \lp  1,j  \rp})\end{eqnarray*}
 by
\[
\ov{z} \mapsto (-\ov{z},\ov{z}).
\]
For ease of notation, we define
\[
\rh_t(1, 1^\prime ):= \bigoplus_{\scriptsize{
\begin{array}{c}  i  \in \lb m \rb \setminus \lp  1  \rp \\  j   \in \lb
n
\rb \setminus \lp 1 \rp \end{array}}}\rh_t(M_{\lb m \rb
\setminus \lp 1, i   \rp,\lb n \rb \setminus \lp  1 , j  \rp}).
\]
and let
\[
\delta:\rh_t(1, 1^\prime ) \rightarrow \rh_t(1^\prime) \oplus \rh_t( 1 )
\] be the unique homomorphism whose restriction to $\rh_t(M_{\lb m
\rb \setminus \lp 1, i  \rp,\lb n \rb \setminus \lp  1 , j
\rp})$ is $\delta^{ i , j  }$ for each pair $( i ,j  )$.

\begin{lem} \label{lescb} The sequence
\beq \\ \nonumber
\cdots \stackrel{\delta}{\rightarrow} \bigoplus_{i \in [m] \setminus \{1\}}
\rh_{t-1}(M_{[m]\setminus\{1,i\},[n]\setminus\{1\}}) \oplus
\bigoplus_{j\in [n] \setminus
\{1\}}\rh_{t-1}(M_{[m]\setminus\{1\},[n]
\setminus\{1,j\}})
\\ \nonumber\stackrel{\phi}{\rightarrow} \,\,\rh_t(M_{m,n}) \,\,\stackrel{\psi} {\rightarrow}
\bigoplus_{\scriptsize \begin{array}{c} i\in [m]\setminus\{1\}\\ j\in[n]\setminus\{1 \}\end{array}}
\rh_{t-2}(M_{[m]\setminus
\{1,i\},[n]\setminus\{1,j\}})
 \stackrel{\delta}{\rightarrow}
\\ \nonumber\bigoplus_{i \in [m]\setminus \{1\}}
\rh_{t-2}(M_{[m]\setminus\{1,i\},[n]\setminus\{1\}}) \oplus
\bigoplus_{j\in [n] \setminus \{1\}}\rh_{t-2}(M_{[m]\setminus\{1\},
[n]\setminus\{1,j\}})
\stackrel{\phi}{\rightarrow}\cdots
\eeq is exact.
\end{lem}

\begin{proof}
Define
\[ X_{m,n}:=\lp G \in M_{m,n}:|(N_G(1) \cup N_G( 1^\prime )) \setminus \lp 1, 1^\prime  \rp|
\leq 1 \rp.
\]
Let $P_{m,n}$ be the subcomplex of $X_{m,n}$ consisting of those $G \in X_{m,n}$ such that
either
$\lp 1,1^\prime \rp \in E(G)$ or both $1$ and $1^\prime$ are isolated in $G$.  As before,
 the  natural projection of chain
complexes induces an isomorphism
\[
\tau:\rh_t(X_{m,n}) \rightarrow \rh_t(X_{m,n},P_{m,n}).
\]

For
$i
\in
\lb m
\rb
\setminus
\lp 1
\rp$ and $j
\in
\lb n
\rb
\setminus
\lp 1
\rp$, let
$$\alpha_{i}: C_t(X_{m,n},P_{m,n}) \to C_{t-1}(M_{[m]\setminus \{1,i\}, [n]
\setminus\{1\}})$$
and $$\alpha^\prime_j: C_t(X_{m,n},P_{m,n}) \to C_{t-1}(M_{[m]\setminus \{1\}, [n]
\setminus\{1,j\}})$$ be the maps defined by letting
$$ \alpha_{i}(x) = \begin{cases} y &\mbox{if } x = 1^\prime i \wedge y
 \mbox{ for some } y \in
C_{t-1}(M_{\lb m \rb
\setminus \lp 1, i   \rp,\lb n \rb \setminus \lp  1   \rp})\\ 0
&\mbox{otherwise,}
\end{cases}$$
and $$ \alpha^\prime_{j}(x) = \begin{cases} y &\mbox{if } x = 1 j^\prime \wedge y
\mbox{ for some } y \in
C_{t-1}(M_{\lb m \rb
\setminus \lp 1  \rp,\lb n \rb \setminus \lp  1, j   \rp})\\ 0
&\mbox{otherwise,}
\end{cases}$$
for each oriented simplex $x$.
 It is straightforward to show that the induced maps,
$$\alpha_i:\rh_t(X_{m,n},P_{m,n}) \to
\rh_{t-1}(M_{[m]\setminus\{1,i\},[n]\setminus\{1\}})$$ and
$$\alpha^\prime_j:\rh_t(X_{m,n},P_{m,n}) \to
\rh_{t-1}(M_{[m]\setminus\{1\},[n]\setminus\{1,j\}})$$ given by
$\alpha_i(\bar z) = \ov{\alpha_i(z)}$ and $\alpha^\prime_j(\bar z) =
\ov{\alpha^\prime_j(z)}$, are well-defined homomorphisms, as is the map
$$\alpha:\rh_t(X_{m,n},P_{m,n}) \rightarrow \rh_{t-1}(1^\prime) \oplus \rh_{t-1}(1)$$
 defined by
\[
\ov{z} \mapsto ((\alpha_{i}(\bar z)),(\alpha^\prime_{j}(\bar z))).
\]
The map $\alpha$ has an inverse analogous to the inverse
$\gamma$ defined in Section \ref{les1}.  Therefore, $\alpha$ is an isomorphism.

For $i\in [m]\setminus\{1\}$, and $j
\in
\lb n
\rb
\setminus
\lp 1
\rp$, the map  $$ \beta_{i,j}:\rh_t(M_{m,n},X_{m,n}) \rightarrow
\rh_{t-2}(M_{[m]\setminus \{1,i\},\lb n
\rb \setminus \lp 1,j \rp})
$$
 induced by the restriction of $\psi_{i,j}$ to $C_t(M_{m,n},X_{m,n})$ is  a
well-defined homomorphism.
 Define
\[
\beta:\rh_t(M_{m,n},X_{m,n}) \rightarrow \bigoplus_{\scriptsize\begin{array}{c} i\in
[m]\setminus\{1\} \\ j \in [n]\setminus\{1\}\end{array}}\!\!\rh_{t-2}(M_{[m]
\setminus\{1,i\},\lb n
\rb \setminus \lp 1,j \rp})
\] by
\[
\ov{z} \mapsto ( \beta_{i,j}(\ov{z})).
\]  As in Section \ref{les1}, $\beta$ is a well-defined isomorphism.

The diagram
\beq \nonumber
& &  \cdots \stackrel{\partial_\ast}{\rightarrow}
\rh_t(X_{m,n}) \quad\stackrel{i_\ast}{\longrightarrow}\quad\rh_t(M_{m,n})
\stackrel{\pi_\ast}{\longrightarrow} \rh_t(M_{m,n},X_{m,n})
\stackrel{\partial_\ast}{\rightarrow} \cdots
\\\nonumber & &
\qquad \qquad\quad\big{\downarrow} \scriptsize{\alpha\tau}
\qquad\qquad\qquad \big{\downarrow} \scriptsize{\mbox{id}}
\qquad\qquad\qquad \qquad \big{\downarrow} \scriptsize{\beta}
\\ \nonumber & & \hspace{-.3in}
\cdots \stackrel{\delta}{\rightarrow} \rh_{t-1}(1) \oplus
\rh_{t-1}( 1^\prime ) \stackrel{\phi}{\rightarrow} \,\,\rh_t(M_{m,n})
\,\,\stackrel{\psi}{\longrightarrow}\, \rh_{t-2}(1, 1^\prime )\quad
\stackrel{\delta}{\longrightarrow} \cdots
\eeq commutes, which yields the result.
\end{proof}

\subsection{The tail end}

For our purposes, we   need only the tail end  of each  long exact sequence.
Recall that $$\nu_n = \lfloor {n+1\over 3} \rfloor -1.$$

\begin{lem} \label{ontom} Let $ \phi$ and $ \psi$ be as in Lemma \ref{lesbouc}.
\begin{enumerate}
\item[(i)] If $n \equiv 0,1 \bmod 3$ then the following is an exact sequence
\beq \nonumber
\bigoplus_{\scriptsize\begin{array}{c} a\in\{1,2\}\\ i\in\{3,\dots,n\}\end{array}}
\rh_{\nu_{n-3}}(M_{[n]\setminus\{1,2,i\}}) \stackrel{\phi}{\rightarrow}
\rh_{\nu_n}(M_n) {\rightarrow} 0.\eeq

\vspace{.2in}
\item[(ii)] If  $n \equiv 2 \bmod 3$ then the following is an exact sequence
\begin{equation*}  \begin{split}
\bigoplus_{\scriptsize\begin{array}{c} a\in\{1,2\}\\
i\in\{3,\dots,n\}\end{array}}\rh_{\nu_{n-3}}(M_{[n]\setminus\{1,2,i\}})
&\stackrel{\phi}{\rightarrow}
\rh_{\nu_n}(M_n) \\ &\stackrel{\psi}{\rightarrow}
\bigoplus_{i\ne j \in \{3,\dots,n\}}\rh_{\nu_{n-4}}(M_{[n]\setminus\{1,2,i,j\}})
\rightarrow 0.\end{split}
\end{equation*}
\end{enumerate}
\end{lem}

\begin{proof} First note that $\nu_{n-3} = \nu_{n}-1$ for all $n$.  Hence  the
  sequence of (i) is a piece of the long exact sequence of  Lemma~\ref{lesbouc},
provided that
$\rh_{\nu_n -2}(M_{n-4}) = 0$.  This follows from (\ref{blvzeq1}), since $\nu_n-2 <
\nu_n-1 =
\nu_{n-4} $ when $n \equiv 0,1 \bmod 3$.

If $n \equiv 2 \bmod 3$, we have that $ \nu_{n-4} = \nu_n-2  $. Hence the sequence of (ii) is
a piece of the  long exact sequence of Lemma~\ref{lesbouc}, by (\ref{blvzeq1}) and the fact
that $\nu_n-2 < \nu_{n-3}$.
\end{proof}

Now recall that, $$\nu_{m,n} = \min\{m,n,\lfloor {m+n+1\over 3} \rfloor\} -1.$$
Note that if $m \le n$ then
\beq \label{note} \nu_{m,n} = \begin{cases} \lfloor {m+n+1\over 3} \rfloor -1 &\mbox{ if } n
\le 2m-1
\\ m-1 & \mbox{ if } n \ge 2m-1, \end{cases} \eeq
and if $n < 2m-1$ then
\beq\label{note2} \nu_{m,n} < m-1. \eeq

\begin{lem}  \label{chessonto} Suppose $m \le n < 2m-1 $. Let $ \phi$ and $ \psi$ be as in
Lemma \ref{les2}.
\begin{enumerate}
\item[(i)] If $m+n \equiv 0,1 \bmod 3$ then
\begin{eqnarray*}
 \bigoplus_{i \in [m]\setminus\{1\}}
\rh_{\nu_{m-2,n-1}}(M_{[m]\setminus\{1,i\},[n]\setminus\{1\}})        &\oplus&
\bigoplus_{j \in
[n]\setminus\{1\}}\rh_{\nu_{m-1,n-2}}(M_{[m]\setminus\{1\},
[n]\setminus\{1,j\}})
\\ \nonumber\stackrel{\phi}{\rightarrow}& & \rh_{\nu_{m,n}}(M_{m,n}) \to 0
\end{eqnarray*} is exact.

\vspace{.2in}
\item[(ii)] If $m+n \equiv 2 \bmod 3$ then
\begin{eqnarray*}
 \bigoplus_{i\in [m]\setminus\{1\}}
\rh_{\nu_{m-2,n-1}}(M_{[m]\setminus\{1,i\},[n]\setminus\{1\}})        &\oplus&
\bigoplus_{j\in
[n]\setminus\{1\}}\rh_{\nu_{m-1,n-2}}(M_{[m]\setminus\{1\},[n]\setminus
\{1,j\}})
\\ \nonumber\stackrel{\phi}{\rightarrow} \rh_{\nu_{m,n}}(M_{m,n})
\stackrel{\psi}{\rightarrow}
\bigoplus_{\scriptsize\begin{array}{c}  i \in [m]\setminus\{1\} \\ j \in
[n]\setminus\{1\}
\end{array}} & &\hspace{-.5in}
\rh_{\nu_{m-2,n-2}}(M_{[m]
\setminus \{1,i\},[n]\setminus\{1,j\}}) \rightarrow 0,
\end{eqnarray*} is exact.
\end{enumerate}
\end{lem}

\begin{proof} Note that for all $m,n$ such that $m \le n < 2m-1$,  \beq \label{nu}
\nu_{m-2,n-1}  =
\nu_{m-1,n-2} =
\nu_{m,n}-1,\eeq
and  \begin{eqnarray*} \nu_{m-2,n-2} = \lfloor{ m+n\over 3}\rfloor-2.
\end{eqnarray*} It follows that if  $m+n \equiv 0,1 \bmod 3$ then
\begin{eqnarray}\label{01} \nu_{m-2,n-2}  &=&  \lfloor{ m+n+1\over 3}\rfloor-2
\\ \nonumber &=&\nu_{m,n}-1.\end{eqnarray} Hence by (\ref{blvzeq2}), we have
$\rh_{\nu_{m,n}-2}(M_{[m]\setminus
\{1,i\},[n]\setminus\{1,j\}}) =0$, which together with (\ref{nu})
implies that the sequence in (i) is a piece of the long exact sequence of Lemma~\ref{lescb}.

If $m+n \equiv 2 \bmod 3$ then
\begin{eqnarray} \nonumber \nu_{m-2,n-2}  &=&  \lfloor{ m+n+1\over 3}\rfloor-3\\ \nonumber&=&\nu_{m,n}-2.
\end{eqnarray} It follows from this, (\ref{blvzeq2}), and (\ref{nu}) that the sequence in
(ii) is a piece of the long exact sequence of Lemma~\ref{lescb}.
\end{proof}

Lemma~\ref{ontom} (resp., \ref{chessonto}) will be used to decompose generators of
$\rh_{\nu_n}(M_n)$ (resp., $\rh_{\nu_{m,n}}(M_{m,n})$) into wedge products of smaller
cycles.  An easy instance of this is given in the next lemma.

\begin{lem}\label{gen01} Suppose $n \equiv 0,1 \bmod 3$.  Then $\rh_{\nu_n}(M_n)$ is
generated by elements of the form \begin{eqnarray*}
(\sigma(1)\sigma(2)-\sigma(1)\sigma(3) ) \,\, &\wedge& \,\,
(\sigma(4)\sigma(5)-\sigma(4)\sigma(6) ) \,\, \wedge \,\, \dots \\ \,\, \dots &\wedge&
\,\,(\sigma(N-2)\sigma(N-1)-\sigma(N-2)\sigma(N) ),
\end{eqnarray*}   where  $\sigma \in \mathfrak S_{n}$ and $N=3 \lfloor {n
\over 3}\rfloor$.
\end{lem}
\begin{proof}  This follows from Lemma~\ref{ontom} (i) by induction on $n$.
\end{proof}

\section{Proof of the BLVZ conjecture} \label{BLVZsec}

  Lemma~\ref{ontom}  is the main tool  of  Bouc's proof of the  conjecture of
Bj\"orner, Lov\'asz, Vr\'ecica and {\v Z}ilvaljevi\'c that $\tilde
H_{\nu_n}(M_n)$ does not vanish.    Bouc first establishes
nonvanishing homology in the most difficult case, the $n \equiv
1\bmod 3$ case.  He then observes that  Lemma~\ref{ontom} enables
one to deduce the $n \equiv 0\bmod 3$ case  from the  $n \equiv
1\bmod 3$ case.  Although not explicitly mentioned by Bouc,  the
same is true for the  remaining $2 \bmod 3$ case.  Indeed,
consider the surjective map $\psi$ of  Lemma~\ref{ontom} (ii).
Since $n-4 \equiv 1 \bmod 3$, the range of $\psi$ does not vanish.
Hence neither does the domain $\tilde H_{\nu_n}(M_n)$.

We now prove the
conjecture for the chessboard complex.

\begin{thm}[Bj\"orner-Lov\'asz-Vr\'ecica-{\v Z}ilvaljevi\'c Conjecture] \label{conj}
\qquad \newline For  $n \ge 3$,
\beq \label{conjeq1} \tilde H_{\nu_n}(M_n) \ne 0, \eeq and for
 $ m+n \ge 3$,
\beq \label{conjeq2} \tilde H_{\nu_{m,n}}(M_{m,n}) \ne 0.
\eeq
\end{thm}

\begin{proof}[Proof of (\ref{conjeq2})]  If $n \ge 2m-1$, then the result follows from
Theorem~\ref{FH}.  So assume that
$m \le n< 2m-1 $.

We will begin with the case that $m+n \equiv 0 \bmod 3$.  The argument for $m+n
\equiv 1 \bmod 3$ is similar and will be left to the reader.   We will use the fact that
$\rh_{\nu_{m+n}}(M_{m+n})$ does not vanish to prove that
$\rh_{\nu_{m,n}}(M_{m,n})$ does not vanish.    Since the chessboard complex $M_{m,n}$ is a
subcomplex of the matching complex
$M_{[m]\uplus [n]^\prime}$,  any cycle
$z$ of
$M_{[m]\uplus [n]^\prime}$ that is  in the chain space of
$M_{m,n}$ must be a cycle in $M_{m,n}$.   Moreover, if $z$ is a boundary in the subcomplex
$M_{m,n}$ then it is also a boundary in
$M_{[m]\uplus [n]^\prime}$.

 Let $k = \frac{2n-m}{3}$. It follows from  $m+n \equiv 0 \bmod 3$, that $k$  is an integer.
The cycle
\begin{eqnarray*} z:= (1\,1^\prime-1\, 2^\prime ) \,\wedge\, (2 \, 3^\prime - 2
\, 4^\prime) \, \wedge & \dots & \wedge \, (k\,(2k-1)^\prime - k\,(2k)^\prime)
\wedge
\\ ((2k+1)^\prime(k+1) -  (2k+1)^\prime (k+2))\wedge  & \dots & \wedge
 \,(n^\prime\, (m-1)-n^\prime\, m )
\end{eqnarray*} of $M_{[m]\uplus[n]^\prime}$ is not a boundary since it is one of the
generators given by Lemma~\ref{gen01}.  Indeed, if any one
of the cycles given by Lemma~\ref{gen01} is a boundary, they all are, which is impossible
since $\rh_{\nu_{m+n}}(M_{[m]\uplus[n]^\prime}) \ne 0$.   The cycle $z$  is  clearly in the
$(\frac{m+n}{ 3} -1)$-chain space of
$M_{m,n}$.  So it is a $(\frac{m+n}{ 3} -1)$-cycle of
$M_{m,n}$ that is not a boundary.  Since by (\ref{note}), $\nu_{m,n} = \nu_{m+n} = \frac{m+n}{
3} -1$, we have  $\rh_{\nu_{m,n}}( M_{m,n}) \ne 0$.

Now suppose $m+n \equiv 2 \bmod 3$.  Just as for the matching complex, the $2 \bmod 3 $ case
is a consequence of the $1 \bmod 3$ case. We use Lemma~\ref{chessonto} (ii).  Since
 $m+n-4 \equiv 1 \bmod 3$,  we have that the range of the surjection $\psi$  does not vanish,
by the $1\bmod 3$ case.  Hence, neither does the domain,
$\rh_{\nu_{m,n}}( M_{m,n})$.
\end{proof}

\section {Torsion in the matching complex}\label{mtorsec}

In this section we  prove Theorem~\ref{mtor}.    We begin with the
following lemma.

\begin{lem}  \label{5-cycle} Suppose $n \equiv 2 \bmod 3$ and $n \ge 5$.  Then
$\rh_{\nu_n}(M_n)$ is generated by elements of the form $\gamma \land \rho$,  where $\gamma
\in \rh_1(M_S)$,  $\rho \in \rh_{\nu_{n-5}}(M_{[n]-S})$, and
$|S| =5$.
\end{lem}

\begin{proof}  The proof is by induction on $n$.  The base step $n =5$ is trivial.  Let
$ n \ge 8$.  For  distinct elements $i,j \in [n]$, recall the map
$$\psi_{i,j}:
\tilde H_{\nu_n}(M_n) \to
\tilde H_{\nu_{n-4}}(M_{[n]\setminus \{1,2,i,j\}})$$ defined in Section~\ref{les1}.  Since
$n-4
\equiv 1 \bmod 3$, it follows from  Lemma~\ref{gen01} that
$$\tilde H_{\nu_{n-4}}(M_{[n]\setminus \{1,2,i,j\}}) = \langle \bar \rho : \rho \in
 Z_{\nu_{n-5}}(M_{[n]\setminus\{1,2,i,j,r\}}),r \in [n]\setminus\{1,2,i,j\}\rangle. $$
 Therefore if $\zeta \in
\rh_{\nu_n}(M_n)$ then
\beq
\label{psi2} \psi_{i,j}(\zeta)=\sum_{r \in [n]\setminus\{1,2,i,j\}} \bar
\rho_{i,j,r},
\eeq
for some $\rho_{i,j,r}\in Z_{\nu_{n-5}}(M_{[n]\setminus\{1,2,i,j,r\}})$.

For distinct elements $a,b,r
\in [n]
\setminus
\{1,2\}$, let
$\gamma_{a,b,r}$ be the cycle
$$(1a \land 2b) +(2b\land ra) + (ra \land 12) + (12 \land rb) + (rb \land 1a)$$  in
$Z_1(M_{\{1,2,a,b,r\}})$.  Clearly
$\gamma_{a,b,r}
\land  \rho_{a,b,r} \in  Z_{\nu_n}(M_n)$ and
\beq
\label{psi1}  \psi_{i,j}(\gamma_{a,b,r}
\land  \rho_{a,b,r}) =\begin{cases} \rho_{i,j,r} &\mbox{if } (i,j)=(a,b) \\ 0
&\mbox{otherwise.} \end{cases}
\eeq
It follows from
(\ref{psi2}) and (\ref{psi1}) that
\begin{eqnarray*}\psi_{i,j}(\zeta  - \hspace{-.1in}\sum_{a\ne b \in [n]\setminus
\{1,2\}}\hspace{-.2in}& &\sum_{r
\in [n]\setminus\{1,2,a,b\}}
\ov{
\gamma_{a,b,r} \land  \rho_{a,b,r}}) \\ &=&
\psi_{i,j}(\zeta) - \sum_{r \in
[n]\setminus\{1,2,i,j\}}
\bar \rho_{i,j,r} \\ &=& 0.\end{eqnarray*}
Hence,
by Lemma~\ref{ontom} (ii), we have
\beq \label{lincomb} \qquad\zeta - \sum_{{a\ne b\in [n]\setminus\{1,2\}}}\,\,\sum_{r \in
[n]\setminus\{1,2,a,b\}} \ov{
\gamma_{a,b,r} \land  \rho_{a,b,r}} \in \ker(\psi)= \im(\phi).\eeq

Clearly  $\im(\phi)$ is generated by elements of the  form $\alpha \land \tau$, where
$\alpha \in \rh_0(M_T)$,  $\tau
\in
\rh_{\nu_{n-3}}(M_{[n]-T})$, and $|T| =3$.  By induction
$\rh_{\nu_{n-3}}(M_{[n]-T})$ is generated by elements
of the form
$\gamma \land \omega$, where
$\gamma \in
\rh_1(M_S)$, $\omega \in \rh_{\nu_{n-8}}(M_{[n]-T-S})$,  and $|S|=5$.  It follows
that  $\im(\phi)$ is generated by elements of the form
$\alpha \land \gamma \land \omega$, where $\alpha \in \rh_0(M_T)$,
$\gamma \in
\rh_1(M_S)$, $\omega \in \rh_{\nu_{n-8}}(M_{[n]-T-S})$, $|T| =3$, and $|S|=5$.  It now
follows from (\ref{lincomb}) that $\zeta$ is an integral  combination of elements of the form
 $\gamma \land \rho$, where $\gamma \in \rh_1(M_S)$,  $\rho
\in
\rh_{\nu_{n-5}}(M_{[n]-S})$ and $|S| =5$.  Since $\zeta$ was arbitrary, $\rh_{\nu_{n}}(M_n)$
is generated by elements of this form.
\end{proof}

We are now ready to prove Theorem~\ref{mtor} which  is restated here.

\vspace{.1in}
\noindent{\bf Theorem \ref{mtor}.}{\em  \,\,For  $n \ge 12$ (except possibly $n=14 $),
$ \tilde H_{\nu_n}(M_n)$ is a nontrivial elementary $3$-group.
\em}
\vspace{.1in}

\begin{proof} By Theorem \ref{Btor}, we need only prove the result for $n \equiv 0,2 \bmod
3$.    We prove the $n \equiv 0 \bmod 3$ case by induction on $n$.
Table~1.1 provides the base step, $$\rh_{\nu_{12}}(M_{12}) =
\Z_3^{56}.$$  The induction step follows from Lemma~\ref{ontom} (i) and Theorem~\ref{conj},
since the homomorphic image of a nontrivial elementary $3$-group  is either
trivial or is  a
nontrivial elementary $3$-group.

Now let  $n \equiv 2 \bmod 3$ and $ n \ge  17$. By Lemma~\ref{5-cycle},
$\rh_{\nu_n}(M_n)$ is generated by elements of the form $\gamma \land \rho$ where $\gamma \in
\rh_1(M_S)$,  $\rho \in \rh_{\nu_{n-5}}(M_{[n]-S})$, and
$|S| =5$.  Since
$n-|S| \ge 12$ and $n-|S| \equiv 0 \bmod 3$, by the
$0
\bmod 3$  case, $$3 (\gamma \land \rho) = \gamma \land 3 \rho = 0.$$  Hence
$\rh_{\nu_n}(M_n)$ has exponent at most 3.  The result now follows from Theorem~\ref{conj}.
\end{proof}

 We conjecture that the result holds for $n=14$ as well.    In principle, one need only check
this on the computer.  However, at the present time the computer, using the software of
\cite{DHSW},
 produces results only up to $n=13$.  We have the following partial result for
$n=14$.

\begin{thm} $\tilde H_{\nu_{14}}(M_{14})$ is a finite group whose
Sylow $3$-subgroup is nontrivial.
\end{thm}

\begin{proof}   By Theorem~\ref{Btor}(i), we have that $\tilde H_{\nu_{14}}(M_{14})$ is
finite.

Let $n = 17$. It follows from Lemma~\ref{5-cycle},  that $\tilde H_{\nu_n}(M_n)$ is generated by
elements of the form $\gamma
\land
\rho$ where
$\gamma \in
\rh_1(M_S)$,
$\rho \in \rh_{\nu_{n-5}}(M_{[n]-S})$, and
$|S| =5$. By Lemma~\ref{gen01}, $\rh_{\nu_{n-5}}(M_{[n]-S})$ is generated by elements of the
form $\alpha \land \omega$ where $\alpha \in \rh_0(M_T)$,
$\rho \in \rh_{\nu_{n-8}}(M_{[n]-S-T})$, and
$|T| =3$.  It follows that $\tilde H_{\nu_n}(M_n)$ is generated by elements of the form
$\alpha \land \tau$ where $\alpha \in \rh_0(M_T)$,
$\tau \in \rh_{\nu_{n-3}}(M_{[n]-T})$, and
$|T| =3$. By (\ref{conjeq1}), at least one of these generators, say $\alpha \land \tau$, is
nonzero.

  We have
$$e(\alpha
\land
\tau) =
\alpha
\land e
\tau = 0,$$ where $e$ is the exponent of  $\rh_{\nu_{14}}(M_{14})$.  Since
$\alpha \land \tau \ne 0$, it follows from Theorem~\ref{mtor} that
$3$ divides $e$, which implies that there is $3$-torsion in  $\tilde H_{\nu_{14}}(M_{14})$.
\end{proof}

\begin{cor} \label{3-prime-m} The  Sylow 3-subgroup of
$\rh_{\nu_{n}}(M_{n})$  is nontrivial for all $n$ such that
$\rh_{\nu_{n}}(M_{n})$ is finite.
\end{cor}

\section{Torsion in the chessboard complex}\label{ctorsec}

In this section we prove  Theorem~\ref{ctor}.  The general idea
is patterned on the proof of the analogous result for the matching complex,
given in the previous section.  However there is a significant  complication.   Just as for the
matching complex, the tail end of the long exact sequence  will be used to decompose generators
into smaller cycles, but this works only if
$n$ is sufficiently close to $m$.   When $n$ is not sufficiently close to $m$, it is
necessary to understand the top homology of the chessboard complex in order to decompose the
generators.  A study of top homology is conducted in
Section~\ref{top}, where an essential decomposition result, Corollary~\ref{topcor}, is
obtained.  This decomposition result  and the tail end of the long exact sequence  will enable
us to prove the key decomposition result:
\vspace{-.1in}\begin{quotation} For all $m \le n \le 2m-2$ except $(m,n) = (4,4)$, the group
$\rh_{\nu_{m,n}}(M_{m,n})$ is  generated by elements of the form
\beq \nonumber (ij^\prime - ik^\prime) \land \rho,\eeq where $i \in [m]$,
$j, k \in [n]$, and $\rho \in \rh_{\nu_{m-1,n-2}}(M_{[m]
\setminus \{i\},[n] \setminus
\{j,k\}})$.
\end{quotation}

We divide the  proof of Theorem~\ref{ctor} into
three cases
 which are handled in  three separate subsections.  An approach to determining torsion for
all finite
$\rh_{\nu_{m,n}}(M_{m,n})$, not covered by Theorem~\ref{ctor}, is discussed in  the final
subsection.

\subsection{The $1 \bmod 3$ case}

For $i,j \in [m]$ and $k,l \in [n]$, let $$\alpha_{i,k^\prime,l^\prime}:=
ik^\prime - i l^\prime \,\,\in \,\,\tilde H_0(M_{\{i\},\{k,l\}}),$$
 and
$$\beta_{i,j,k^\prime}:= ik^\prime - j k^\prime\,\,\in \,\,
\rh_0(M_{\{i,j\},\{k\}}).$$  We refer to the fundamental cycle $\alpha_{i,k^\prime,l^\prime}$ as
an
$\alpha$-cycle and the fundamental cycle $\beta_{i,j, k^\prime}$ as
a
$\beta$-cycle.  We   also need to  view these fundamental cycles as elements of
$\rh_0(M_{\{i,j\},\{k,l\}})$.

\begin{lem} \label{rel} In  $\rh_0(M_{\{i,j\},\{k,l\}})$ we have
$$\alpha_{ j,k^\prime,l^\prime} = - \alpha_{i,k^\prime,l^\prime} = -\beta_{i,j,k^\prime} =
\beta_{i,j,l^\prime}.$$
\end{lem}

\begin{proof}  The first equation follows from $$\partial((ik^\prime \land jl^\prime) -
(il^\prime \land jk^\prime)) = (ik^\prime - il^\prime) +(jk^\prime -j l^\prime).$$  The
second equation follows from
$$\partial(il^\prime \land jk^\prime) = (il^\prime - ik^\prime) +(ik^\prime -j k^\prime).$$
The  third equation follows from
$$\partial((ik^\prime \land jl^\prime) + (il^\prime \land jk^\prime)) = (ik^\prime -
jk^\prime) +(il^\prime -j l^\prime).$$
\end{proof}

\begin{lem} \label{gener1} Suppose $m+n \equiv 1 \bmod 3$ and $ m \le n \le 2m-2 $. Then
$\rh_{\nu_{m,n}}(M_{m,n})$ is  generated by elements of the form
\beq \label{alpha} \alpha_{i,j^\prime,k^\prime} \land \rho,\eeq where $i \in [m]$,
$j, k \in [n]$, and $\rho \in \rh_{\nu_{m-1,n-2}}(M_{[m]
\setminus \{i\},[n] \setminus
\{j,k\}})$.
\end{lem}

\begin{proof}  First  note that it follows from  Lemma~\ref{chessonto} (i) that
$\rh_{\nu_{m,n}}(M_{m,n})$ is generated by elements of the form given in (\ref{alpha}) and
elements of the form
\beq \label{beta} \beta_{i,j,k^\prime} \land \rho,\eeq where $i,j \in [m]$, $ k \in
[n]$, and $\rho \in \rh_{\nu_{m-2,n-1}}(M_{[m] \setminus
\{i,j\},[n] \setminus
\{k\}})$.

We will show by induction on $m$ that the elements  of the form given in
(\ref{beta}) can be expressed as integral combinations of elements of the form given in
(\ref{alpha}). The base step, $m=n=2$, follows from
Lemma~\ref{rel}.
Now suppose $m >2$.

{\bf Case 1.}   Say $n<2m-2$. Then $n-1 \leq 2(m-2)-2$ and we
apply the induction hypothesis to $\rh_{\nu_{m-2,n-1}}(M_{[m]
\setminus \{i,j\},[n] \setminus \{k\}})$.  By replacing $\rho$ in
(\ref{beta}) by an integral combination of wedge products each of
which contains an $\alpha$-cycle, we are able to express
$\beta_{i,j,k^\prime} \land \rho$ as an integral combination of
wedge products each of which contains an $\alpha$-cycle.

{\bf Case 2.}   Say $n=2m-2$.  Then $n-1>2(m-2)-2$ and it follows
from (\ref{note}) that
$$\nu_{m-2,n-1} = m-3,$$  so we can apply Corollary~\ref{topcor} to
$\rh_{\nu_{m-2,n-1}}(M_{[m] \setminus \{i,j\},[n]
\setminus
\{k\}})$, which implies that generators given in (\ref{beta}) can be expressed as integral
combinations of elements of the form
\beq \label{fundwedge} \rho_{U,V} \land \gamma,\eeq where $|U| = |V|-1$, $\rho_{U,V}\in
\rh_{|U|-1}(M_{U,V})$, and
$\gamma
\in \rh_{\nu_{m,n}- |U|}(M_{[m]\setminus U,[n]\setminus V})$.

We will show that if $|U|>1$ then  \beq
\label{claim}\rh_{\nu_{m,n}-|U|}(M_{[m]\setminus U,[n]
\setminus V})= 0,\eeq
from which it follows that the wedge product in (\ref{fundwedge}) is $0$.  From this
it follows that the generators in given in (\ref{beta}) can be expressed as integral
combinations of generators given  in (\ref{alpha}), since $\rho_{U,V}$ is an $\alpha$-cycle when
$|U|=1$.

Since $n=2m-2$ and $m>2$, we have $n>m$.  Thus $n-|V| \geq m-|U|$.
It follows that
\[
\nu_{m-|U|,n-|V|}=\min(m-|U|,\lfloor\frac{m-|U|+n-|V|+1}{3}\rfloor)-1.
\]

Suppose $|U| >1$.   We will use (\ref{blvzeq2}) of Theorem~\ref{conn} to prove (\ref{claim}).
From (\ref{note2}) we have
\beq
\label{observe} \nu_{m,n} -|U| < m- |U| -1.\eeq

We also need to check that \beq \label{finally}\nu_{m,n} - |U| < \lfloor {m - |U| + n- |V| +1
\over 3} \rfloor -1.\eeq  By (\ref{note}), we have
$$\nu_{m,n} - |U| = {m+n-1  \over 3} -1- |U|.$$ The right side of (\ref{finally}) equals
$$\lfloor {m - |U| + n- |U| -1+1 \over 3} \rfloor -1 = {m+n-1
\over 3} + \lfloor {-2 |U| +1 \over 3}\rfloor -1.$$  So (\ref{finally}) is equivalent to
$$-|U| <  \lfloor {-2 |U| +1 \over 3}\rfloor,$$  which clearly holds when $|U| \ge 2$.
Hence by (\ref{blvzeq2}), equation (\ref{claim}) holds.
\end{proof}

\begin{lem} \label{gen1} Suppose $m+n \equiv 1 \bmod 3$ and $ m \le n \le 2m-2 $. Then
$\rh_{\nu_{m,n}}(M_{m,n})$ is  generated by elements of the form
\begin{eqnarray}\label{eqgener} \alpha_{\sigma(1),\tau(1)^\prime,\tau(2)^\prime}
\land \alpha_{\sigma(2),\tau(3)^\prime,\tau(4)^\prime}
\land
\!\!\!&\cdots&\!\!\!
\land
\alpha_{\sigma(t),\tau(2t-1)^\prime,\tau(2t)^\prime} \land \\ \nonumber
\beta_{\sigma(t+1),\sigma(t+2),\tau(2t+1)^\prime} \land
\beta_{\sigma(t+3),\sigma(t+4),\tau(2t+2)^\prime}
\,\land
\!\!\!&\cdots&\!\!\!
\land \,
\beta_{\sigma(m-2),\sigma(m-1),\tau(n)^\prime}\, ,
\end{eqnarray} where  $\sigma \in \mathfrak S_m$, $\tau \in \mathfrak S_n$, and
$t=\frac{2n-m+1} 3$.
\end{lem}
\begin{proof} We use induction on $m$. When $m=2$, the result is immediate from
Lemma~\ref{gener1}.  When $2 < m < n$ the result follows from Lemma~\ref{gener1} and the
induction hypothesis applied to $\rh_{\nu_{m-1,n-2}}(M_{[m] \setminus
\{i\},[n] \setminus
\{j,k\}})$.

When  $2 < m =n$, we also use Lemma~\ref{gener1} and apply the induction hypothesis to
$\rh_{\nu_{m-1,n-2}}(M_{[m] \setminus \{i\},[n] \setminus
\{j,k\}})$.  However there is an additional step. Since $m +n \equiv 1 \bmod 3$,
we have $5 \le m =n$.  Hence  $ n-2 \le m-1
\le 2(n-2)-2.$  This allows us to apply the induction hypothesis with the role of the
$\alpha$-cycles and the $\beta$-cycles switched.  Hence we have that
$\rh_{\nu_{m,n}}(M_{m,n})$ is generated by elements of the form
\begin{eqnarray*} \alpha_{\sigma(1),\tau(1)^\prime,\tau(2)^\prime} \land
\alpha_{\sigma(2),\tau(3)^\prime,\tau(4)^\prime}
\land \dots
\land
\alpha_{\sigma(t),\tau(2t-1)^\prime,\tau(2t)^\prime} \,\,\,\land \\
\beta_{\sigma(t+1),\sigma(t+2),\tau(2t+1)^\prime} \land
\beta_{\sigma(t+3),\sigma(t+4),\tau(2t+2)^\prime}
\,\land \dots
\land \,
\beta_{\sigma(m-1),\sigma(m),\tau(n-1)^\prime}\, ,
\end{eqnarray*} where $\sigma \in \mathfrak S_m$, $\tau \in
\mathfrak S_n$, and $t =
\frac{2n-m+1} 3-1$.  To complete the proof, we use Lemma~\ref{rel} to change one of the
$\beta$-cycles to an $\alpha$-cycle.
\end{proof}

\begin{thm} \label{chess1} Suppose $m+n \equiv 1 \bmod 3$ and $ m \le n \le 2m-5 $. Then
$\rh_{\nu_{m,n}}(M_{m,n})$ is a cyclic group of order $3$ generated by
\begin{eqnarray}\label{genchess1} \alpha_{1,1^\prime,2^\prime} \land
\alpha_{2,3^\prime,4^\prime} \land \cdots \land
\alpha_{t,(2t-1)^\prime,(2t)^\prime} &\land& \\ \nonumber
\beta_{t+1,t+2,(2t+1)^\prime} \land \beta_{t+3,t+4,(2t+2)^\prime}
\,\land
\!\!\!&\cdots&\!\!\!
\land \,
\beta_{m-2,m-1,n^\prime}\, ,
\end{eqnarray} where  $t =
\frac{2n-m+1} 3$.
\end{thm}

\begin{proof} We  use the relations of Lemma~\ref{rel} to show that the generators of
Lemma~\ref{gen1} are all equal up to sign.   It suffices to  show
that \beq \label{sigma} \\ \nonumber
\alpha_{1,1^\prime,2^\prime}\land \!\!\!&\cdots&\!\!\!  \land
\alpha_{t,(2t-1)^\prime,(2t)^\prime} \land
\beta_{t+1,t+2,(2t+1)^\prime}\land \dots \land
\beta_{m-2,m-1,n^\prime}\\ \nonumber  &= & \sgn(\sigma)\,\,
\alpha_{\sigma(1),1^\prime,2^\prime}\land \dots \land
\alpha_{\sigma(t),(2t-1)^\prime,(2t)^\prime} \\ \nonumber &
&\qquad \qquad\land\,
\beta_{\sigma(t+1),\sigma(t+2),(2t+1)^\prime}\land \dots \land
\beta_{\sigma(m-2),\sigma(m-1),n^\prime}, \eeq and \beq
\label{tau} \\ \nonumber \alpha_{1,1^\prime,2^\prime}\land
\!\!\!&\cdots&\!\!\!  \land \alpha_{t,(2t-1)^\prime,(2t)^\prime}
\land \beta_{t+1,t+2,(2t+1)^\prime}\land \dots \land
\beta_{m-2,m-1,n^\prime}\\ \nonumber  &= & \sgn(\tau)\,\,
\alpha_{1,\tau(1)^\prime,\tau(2)^\prime}\land \dots \land
\alpha_{t,\tau(2t-1)^\prime,\tau(2t)^\prime} \\ \nonumber &
&\qquad \qquad\land\, \beta_{t+1,t+2,\tau(2t+1)^\prime}\land \dots
\land \beta_{m-2,m-1,\tau(n)^\prime}, \eeq for all $\sigma \in
\mathfrak S_m$ and $\tau \in \mathfrak S_n$.

For the sake of ease of notation and getting to the heart of the argument, we prove
(\ref{sigma}) and (\ref{tau}) for  $m = n=5$.  The general argument is essentially the
same.   To prove
\beq \label{sigma5}\\ \nonumber
\alpha_{1,1^\prime,2^\prime}\land \alpha_{2,3^\prime,4^\prime} \land
\beta_{3,4,5^\prime}= \sgn(\sigma) \,\,
\alpha_{\sigma(1),1^\prime,2^\prime}
\land
\alpha_{\sigma(2),3^\prime,4^\prime} \land \beta_{\sigma(3),\sigma(4),5^\prime} \eeq
 for all $\sigma$, it suffices to prove this for $\sigma$ in the  set of transpositions $\{
(1,5), (2,5),
(1,3), (1,4)\}$, which generates $\mathfrak S_5$.

{\bf Case 1.}  $\sigma = (1,5)$.  By Lemma~\ref{rel}, $\alpha_{1,1^\prime,2^\prime} = -
\alpha_{5,1^\prime,2^\prime}$.  Hence
$$\alpha_{1,1^\prime,2^\prime}\land \alpha_{2,3^\prime,4^\prime} \land
\beta_{3,4,5^\prime} = -
\alpha_{5,1^\prime,2^\prime}\land \alpha_{2,3^\prime,4^\prime} \land
\beta_{3,4,5^\prime.}$$

{\bf Case 2.}  $\sigma = (2,5)$.  This is similar to Case 1.

{\bf Case 3.} $\sigma = (1,3)$.  By repeated applications of Lemma~\ref{rel}, we have
\begin{eqnarray*} \alpha_{1,1^\prime,2^\prime}\land \alpha_{2,3^\prime,4^\prime}
\land \beta_{3,4,5^\prime}  &=&
\alpha_{1,1^\prime,2^\prime}\land \beta_{2,5,3^\prime} \land \beta_{3,4,5^\prime} \\ &=&
\alpha_{1,1^\prime,2^\prime}\land \beta_{2,5,3^\prime} \land
\alpha_{4,4^\prime,5^\prime} \\ &=& -\alpha_{3,1^\prime,2^\prime}\land
\beta_{2,5,3^\prime} \land \alpha_{4,4^\prime,5^\prime} \\ &=&
-\alpha_{3,1^\prime,2^\prime}\land \beta_{2,5,3^\prime} \land \beta_{1,4,5^\prime}
\\ &=& -\alpha_{3,1^\prime,2^\prime}\land \alpha_{2,3^\prime, 4^\prime} \land
\beta_{1,4,5^\prime}
\end{eqnarray*}

{\bf Case 4.} $\sigma = (1,4)$.  This is similar to Case 3.

To show\beq \label{switch} \\ \nonumber
\alpha_{1,1^\prime,2^\prime}\land \alpha_{2,3^\prime,4^\prime} \land
\beta_{3,4,5^\prime}= \sgn(\tau) \,\,
\alpha_{1,\tau(1)^\prime,\tau(2)^\prime}
\land
\alpha_{2,\tau(3)^\prime,\tau(4)^\prime} \land \beta_{3,4,\tau(5)^\prime}\eeq
 for all
$\tau \in \mathfrak S_5$, we use Lemma~\ref{rel} to exchange an
$\alpha$-cycle for a $\beta$-cycle. That is,  by Lemma~\ref{rel}, equation (\ref{switch}) is
equivalent to
$$
\alpha_{1,1^\prime,2^\prime}
\land
\beta_{2,5,3^\prime} \land \beta_{3,4,5^\prime} = \sgn(\tau) \,\,
\alpha_{1,\tau(1)^\prime,\tau(2)^\prime}
\land
\beta_{2,\tau(5)^\prime,\tau(3)^\prime} \land \beta_{3,4,\tau(5)^\prime}.$$  This is
equivalent to
 (\ref{sigma5}) with the role of the $\alpha$-cycles and
$\beta$-cycles switched.

It is straightforward to extend the argument for $m = n =5$ to general $m \le n
\le 2m-5$ since $\mathfrak S_m$ is generated by the set of transpositions $\{(1,m)\dots
(t,m), (1,t+1), \dots, (1,m-1)\}$, and  the expressions on each side of (\ref{sigma}) and
(\ref{tau})  contain at least two
$\alpha$-cycles and at least one
$\beta$-cycle.

We  now show that the order of the cyclic group $\rh_{\nu_{m,n}}(M_{m,n})
$ is 3 by induction on $m$.  The base step $\rh_{\nu_{5,5}}(M_{5,5}) = \Z_3$ is given in
Table~1.2.  Let
$m
\ge 6$.  The generator given in (\ref{genchess1}) can  be expressed as
$$\alpha_{1,1^\prime,2^\prime} \land \rho$$ where $\rho \in
\rh_{\nu_{m,n}-1}(M_{[m]\setminus \{1\},[n] \setminus
\{1,2\}})
$.   If $m < n$   then clearly $ m-1 \le n-2 \le 2(m-1) -5$. If  $m = n$ then
$m=n\ge 8$ which implies
$ n-2
\le m-1 \le 2(n-2) -5$.   In either case, $\nu_{m,n}-1 = \nu_{m-1,n-2}$, and we can apply the
induction hypothesis to
$  \rh_{\nu_{m,n}-1}(M_{[m]\setminus \{1\},[n]
\setminus
\{1,2\}})
$  to obtain
$$3 (\alpha_{1,1^\prime,2^\prime} \land \rho ) = \alpha_{1,1^\prime,2^\prime} \land 3 \rho =
0.$$ Since, by Theorem~\ref{conj}, $\rh_{\nu_{m,n}}(M_{m,n})$ is nonvanishing, it has order $3$.
\end{proof}

\subsection{The $0 \bmod 3$ case}
\begin{lem} \label{gener0}Suppose  $m+n \equiv 0 \bmod 3$ and $m \le n \le 2m-3 $. Then
$\rh_{\nu_{m,n}}(M_{m,n})$ is  generated by elements of the form
\beq \label{alpha0} \alpha_{i,j^\prime,k^\prime} \land \rho,\eeq where $i \in [m]$,
$j, k \in [n]$, and $\rho \in \rh_{\nu_{m-1,n-2}}(M_{[m]
\setminus \{i\},[n] \setminus\{j,k\}})$.
\end{lem}

\begin{proof}The proof, although similar to  the proof of  Lemma~\ref{gener1}, requires an
additional step.  By Lemma~\ref{chessonto} (i), we have that
$\rh_{\nu_{m,n}}(M_{m,n})$ is  generated by elements of the form given in (\ref{alpha0}) and
elements of the form
\beq \label{beta0} \beta_{i,j,k^\prime} \land \rho,\eeq where $i,j \in [m]$, $ k \in
[n]$, and $\rho \in \rh_{\nu_{m-2,n-1}}(M_{[m] \setminus
\{i,j\},[n] \setminus
\{k\}})$.  It follows from this that
$\rh_{\nu_{3,3}}(M_{3,3})$ is generated by elements of the form
$\alpha_{i_1,j_1^\prime,j_2^\prime} \land \beta_{i_2,i_3,j_3^\prime} $, which takes care of
the base step of an induction proof.  Now assume $m >3$.

{\bf Case 1.} Say $n<2m-3$.  Then $n-1 \le 2(m-2) -3$.  By
applying the induction hypothesis to
$\rh_{\nu_{m-2,n-1}}(M_{[m]\setminus\{i,j\},[n]\setminus\{k\}})$,
we have that the generators  given in (\ref{beta0}) can be
expressed as  integral combinations of generators given in
(\ref{alpha0}).

{\bf Case 2.} Say $n=2m-3$.  Then $n-1 > 2(m-2) -1$, so by
(\ref{note2}), we have $\nu_{m-2,n-1} = m-3$. By applying
Corollary~\ref{topcor} to
$\rh_{\nu_{m-2,n-1}}(M_{[m]\setminus\{i,j\},[n]\setminus\{k\}})$,
we see that generators given in (\ref{beta0}) can be expressed as
integral combinations of elements of the form \beq \label{fund0}
\rho_{U,V} \land \gamma,\eeq where $|U| = |V|-1$, $\rho_{U,V} \in
\rh_{|U|-1}(M_{U,V})$, and $\gamma \in \rh_{\nu_{m,n}-
|U|}(M_{[m]\setminus U,[n]\setminus V})$.

One can show that if $|U| > 2$ then $\rh_{\nu_{m,n}- |U|}(M_{[m]\setminus U,[n]\setminus V})
= 0$ by using an argument similar to the one that was used to prove (\ref{claim}).  We leave
the straightforward details to the reader.  This allows us to conclude that
$\rh_{\nu_{m,n}}(M_{m,n})$ is  generated by elements given in (\ref{alpha0}) and
(\ref{fund0}), where $2=|U| = |V|-1$.

We now show that any generator of the form given in (\ref{fund0}),
where $(|U|,|V|) = (2,3)$, can be expressed as  integral
combination of generators given in (\ref{alpha0}), which will
complete the proof.  Since $m>3$ and $n=2m-3$, we have $m<n$. Thus
$$m-2 \le n-3 \le 2(m-2) -2.$$ By (\ref{note}), we have $\nu_{m,n}
- |U| = \nu_{m-|U|,n-|V|}$. It therefore follows from
Lemma~\ref{gener1}, that $\rh_{\nu_{m,n}- |U|}(M_{[m]\setminus
U,[n]\setminus V})$ is generated by wedge products that contain an
$\alpha$-cycle.
\end{proof}

The next result follows easily from Lemma~\ref{gener0} by induction.

\begin{lem} \label{gen0} Suppose  $m+n \equiv 0 \bmod 3$ and $m \le n \le 2m-3 $. Then
$\rh_{\nu_{m,n}}(M_{m,n})$ is  generated by elements  of the form \beq
\label{type1} \alpha_{i_1,j_1^\prime,j_2^\prime} \land \beta_{i_2,i_3,j_3^\prime}
\land \xi,\eeq where $i_1,i_2,i_3 \in [m]$, $j_1,j_2,j_3 \in [n]$
and
$$\xi \in \rh_{\nu_{m,n}-2}(M_{[m] \setminus\{i_1,i_2,i_3\},[n]
\setminus\{j_1,j_2,j_3\}}).$$
\end{lem}

\vspace{.2in}
For distinct $i_1,i_2,i_3 \in [m]$ and distinct $j_1,j_2,j_3
\in [n]$, let
$$u_{i_1,i_2,j_1^\prime,j_2^\prime,j_3^\prime} := $$ $$ i_1j^\prime_1
\land i_2j_2^\prime + i_2j_2^\prime \land i_1 j_3^\prime + i_1 j_3^\prime \land i_2
j_1^\prime + i_2 j_1^\prime\land i_1 j_2^\prime + i_1 j_2^\prime \land i_2 j_3^\prime + i_2
j_3^\prime \land i_1 j_1^\prime$$ and
$$v_{i_1,i_2,i_3,j_1^\prime,j_2^\prime}:=$$
$$i_1j^\prime_1
\land i_2j_2^\prime + i_2j_2^\prime \land i_3 j_1^\prime + i_3 j_1^\prime \land i_1
j_2^\prime + i_1 j_2^\prime\land i_2 j_1^\prime + i_2 j_1^\prime \land i_3 j_2^\prime + i_3
j_2^\prime \land i_1 j_1^\prime.$$
We shall view $u_{i_1,i_2,j_1^\prime,j_2^\prime,j_3^\prime}$ and
$v_{i_1,i_2,i_3,j_1^\prime,j_2^\prime}$ as  elements of
$\rh_{\nu_{3,3}} (M_{\{i_1,i_2,i_3\},\{j_i,j_2,j_3\}})$ as well as of $\rh_{\nu_{2,3}}
(M_{\{i_1,i_2\},\{j_i,j_2,j_3\}}) $  and  $\rh_{\nu_{3,2}}
(M_{\{i_1,i_2, i_3\},\{j_i,j_2\}}) $, respectively.

\begin{lem}\label{rel0} In $\rh_1(M_{3,3})$ we have,
 $$ 3( \alpha_{1,1^\prime,2^\prime}\land\beta_{2,3,3^\prime}) =
- u_{2,3,1^\prime,2^\prime,3^\prime}-
 v_{1,2,3,1^\prime,2^\prime}  - 2 ( v_{1,2,3,2^\prime,3^\prime} +
 u_{1,2,1^\prime,2^\prime,3^\prime} ).$$
\end{lem}

\begin{proof} It is straightforward to verify that
\begin{eqnarray*} \partial ( 11^\prime \land 22^\prime \land 33^\prime + 12^\prime
\land 23^\prime \land 31^\prime + 12^\prime \land 21^\prime \land 33^\prime + 11^\prime \land
32^\prime \land 23^\prime) \\ =
u_{2,3,1^\prime,2^\prime,3^\prime} +
v_{1,2,3,1^\prime,2^\prime} - \alpha_{1,1^\prime,2^\prime}
\land \beta_{2,3,3^\prime} - 2 (\alpha_{3,2^\prime,3^\prime} \land
\beta_{1,2,1^\prime}).
\end{eqnarray*}
  Consequently, in
$\rh_1(M_{3,3})$,$$\alpha_{1,1^\prime,2^\prime}
\land \beta_{2,3,3^\prime} =   u_{2,3,1^\prime,2^\prime,3^\prime} +
 v_{1,2,3,1^\prime,2^\prime}  - 2 (\alpha_{3,2^\prime,3^\prime} \land
\beta_{1,2,1^\prime}).
 $$ By symmetry (exchanging $\alpha$ with $\beta$, $u$ with $v$, and $i$ with $i^\prime$),
$$\beta_{1,2,1^\prime} \land \alpha_{3,2^\prime,3^\prime} =
 v_{1,2,3,2^\prime,3^\prime} +   u_{1,2,1^\prime,2^\prime,3^\prime} -
2(\beta_{2,3,3^\prime} \land \alpha_{1,1^\prime,2^\prime}).$$ By substituting the second
equation into the first equation, we get
\begin{eqnarray*}\alpha_{1,1^\prime,2^\prime}\!\! &\land& \!\!\beta_{2,3,3^\prime} = \\
 u_{2,3,1^\prime,2^\prime,3^\prime}\!\!&  + &\!\!
 v_{1,2,3,1^\prime,2^\prime}  + 2 ( v_{1,2,3,2^\prime,3^\prime} +
 u_{1,2,1^\prime,2^\prime,3^\prime} - 2(\beta_{2,3,3^\prime} \land
\alpha_{1,1^\prime,2^\prime})),
\end{eqnarray*} which implies that
$$ 3( \alpha_{1,1^\prime,2^\prime}\land\beta_{2,3,3^\prime}) =
- u_{2,3,1^\prime,2^\prime,3^\prime}-
 v_{1,2,3,1^\prime,2^\prime}  - 2 ( v_{1,2,3,2^\prime,3^\prime} +
 u_{1,2,1^\prime,2^\prime,3^\prime} ).$$
\end{proof}

\begin{thm} \label{chess0} Suppose $m+n \equiv 0 \bmod 3$ and $ m \le n \le 2m-9 $.
Then $\rh_{\nu_{m,n}}(M_{m,n})$ is a nontrivial
$3$-group of exponent at most 9.
\end{thm}

\begin{proof} It follows from Lemmas~\ref{gen0} and \ref{rel0}, that $3
\rh_{\nu_{m,n}}(M_{m,n})$ is generated by elements of the form $$\rho_{U,V} \land
\omega,$$ where $( |U|, |V|) \in \{ (2,3),(3,2)\}$, $\rho_{U,V} \in
\rh_1(M_{U,V})$, and
 $\omega
\in
\rh_{\nu_{m,n}- 2}(M_{[m]\setminus U,[n]\setminus V})$.  We can show that
\beq
\label{exp3} 3(\rho_{U,V} \land \omega) =\rho_{U,V} \land 3 \omega = 0\eeq by applying
Theorem~\ref{chess1}, if we first check that
$m-|U|$ and $n-|V|$ satisfy the hypothesis of the theorem.  Clearly $m-|U|+n-|V| = m+n-5
\equiv 1 \bmod 3$.  We leave it to the reader to check the inequalities in each of the three
cases:
\begin{enumerate} \item  $ m<n$ and $( |U|, |V|) = (2,3)$
\item $m=n$ and $( |U|, |V|) =
(2,3)$
\item $m \le n$ and
$( |U|, |V|) = (3,2)$.
\end{enumerate} It  follows from (\ref{exp3}) that $\rh_{\nu_{m,n}}(M_{m,n})$
has exponent at most $9$, and from Theorem~\ref{conj} that the group is nontrivial.\end{proof}

\subsection{The $2 \bmod3$ case}

\begin{lem} \label{alphabeta} Suppose  $m+n \equiv 2 \bmod 3$ and $4\le m \le n \le 2m-4 $.
Then  $\rh_{\nu_{m,n}}(M_{m,n})$ is  generated by elements of the form
\beq  \alpha_{i,j^\prime,k^\prime} \land \rho,\eeq where $i \in [m]$, $j, k \in
[n]$, and $\rho \in \rh_{\nu_{m-1,n-2}}(M_{[m] \setminus
\{i\},[n] \setminus
\{j,k\}})$, and elements of the form
\beq \label{beta2} \beta_{i,j,k^\prime} \land \rho,\eeq where $i,j \in [m]$, $ k \in
[n]$, and $\rho \in \rh_{\nu_{m-2,n-1}}(M_{[m] \setminus
\{i,j\},[n] \setminus
\{k\}})$.
\end{lem}

\begin{proof} We claim  that $\bigoplus_{i ,j }
\rh_{\nu_{m-2,n-2}}(M_{[m]
\setminus \{1,i\},[n]\setminus\{1,j\}})$ is generated by elements of the
form $\psi(\alpha_{r,s^\prime,t^\prime}
\land
\rho ),$ where $\psi$ is the surjection of Lemma~\ref{chessonto}~(ii), and
\begin{itemize}
\item $ r \in [m]\setminus \{1\}$
\item  $s, t  \in [n]\setminus \{1\}$
\item  $\rho \in
\rh_{\nu_{m-1,n-2}}(M_{[m]
\setminus\{r\},[n]\setminus
\{s,t\}})$.
\end{itemize}
We prove this claim by first using Lemma~\ref{gener1} to observe that
$$
\rh_{\nu_{m-2,n-2}}(M_{[m]
\setminus \{1,i\},[n]\setminus\{1,j\}})$$ is generated by elements of the
form $\alpha_{r,s^\prime,t^\prime} \land
\tau,$ where
\begin{itemize}
\item  $r\in [m]
\setminus \{1,i\}$
\item  $s,t \in [n] \setminus
\{1,j\}$
\item $\tau \in \rh_{\nu_{m-3,n-4}}(M_{[m]\setminus \{1,i,r\},[n]\setminus
\{1,j,s,t\}})$.
\end{itemize}   The map
$$\psi:\rh_{\nu_{m-1,n-2}}(M_{[m] \setminus\{r\},[n]\setminus
\{s,t\}}) \to
\bigoplus_{i,j} \rh_{\nu_{m-3,n-4}}(M_{[m]\setminus
\{1,i,r\},[n]\setminus \{1,j,s,t\}})$$ is surjective by
Lemma~\ref{chessonto} (ii).  Hence for  $$\tau \in
\rh_{\nu_{m-3,n-4}}(M_{[m]\setminus
\{1,i,r\},[n]\setminus
\{1,j,s,t\}}),$$ we can let $\rho \in
\rh_{\nu_{m-1,n-2}}(M_{[m]
\setminus\{r\},[n]\setminus
\{s,t\}})$  be such that $\psi(\rho) = \tau.$  It follows directly from the
definition of $\psi$ that
$$\psi(\alpha_{r,s^\prime,t^\prime} \land
\rho) = \alpha_{r,s^\prime,t^\prime} \land
\tau,$$ which proves our claim.

 Let
$\gamma
\in
\rh_{\nu_{m,n}}(M_{m,n})$. We express $\psi(\gamma)$ as an integral combination of generators:
\beq \nonumber\psi(\gamma) = \sum_{r,s,t,\rho} c_{r,s,t,\rho}
\psi(\alpha_{r,s^\prime,t^\prime} \land \rho) =
\psi\left(
\sum_{r,s,t,\rho} c_{r,s,t,\rho} (\alpha_{r,s^\prime,t^\prime} \land \rho)\right),
\eeq for some $c_{r,s,t,\rho} \in \Z$.  It follows from Lemma~\ref{chessonto} (ii)  that
$$\gamma - \sum_{r,s,t,\rho} c_{r,s,t,\rho} (\alpha_{r,s^\prime,t^\prime} \land
\rho) \in \im  \phi.$$   Hence
$\gamma$ can be expressed as an integral combination  of  elements of the form given in the
statement of the lemma.
\end{proof}

Next we show that the elements given in (\ref{beta2}) can be removed from the generating set.

\begin{lem} \label{gener2} Suppose  $m+n \equiv 2 \bmod 3$ and $5\le m \le n \le 2m-4 $.
Then  $\rh_{\nu_{m,n}}(M_{m,n})$ is  generated by elements of the form
\beq \label{alpha2} \alpha_{i,j^\prime,k^\prime} \land \rho,\eeq where $i \in [m]$,
$j, k \in [n]$, and $\rho \in \rh_{\nu_{m-1,n-2}}(M_{[m]
\setminus \{i\},[n] \setminus
\{j,k\}})$.
\end{lem}

\begin{proof} The proof is similar to  the proofs of  Lemmas~\ref{gener1} and
\ref{gener0}.  We use induction on $m$.  The base step,
$(m,n) = (5,6)$, is part of Case 2 below, which does not require the induction hypothesis.

We will  show that generators  given in (\ref{beta2})
can be expressed as integral combinations of generators  given in (\ref{alpha2}).

{\bf Case 1.} Say $n<2m-4$.  Then $5 <m$ and  $ n-1 \le 2(m-2)
-4$. Moreover, $m \ne 6$ because otherwise $n-1 \le 4$.  Hence,
$5 \le m-2 \le n-1 \le 2(m-2) -4$, which enables us to apply the
induction hypothesis to $\rh_{\nu_{m-2,n-1}}(M_{[m] \setminus
\{i,j\},[n] \setminus \{k\}})$.

{\bf Case 2.} Say $n=2m-4$. Since $m \ge 5$, it follows that $n
> m$.

 By
Corollary~\ref{topcor} applied to
$\rh_{\nu_{m-2,n-1}}(M_{[m]
\setminus \{i,j\},[n] \setminus
\{k\}})$,  the generators given in (\ref{beta2}) can be expressed as integral
combinations of elements of the form
\beq \label{fundwedge2} \rho_{U,V} \land \gamma,\eeq where $|U| = |V|-1$, $\rho_{U,V} \in
\rh_{|U|-1}(M_{U,V})$, and
$\gamma \in \rh_{\nu_{m,n}- |U|}(M_{[m]\setminus U,[n]\setminus V})$.

An argument similar to the one used in the proof of
Lemma~\ref{gener1} shows that if $|U|>4$, then the wedge product
in (\ref{fundwedge2}) is $0$.  From this it follows that the
generators in given in (\ref{beta2}) can be expressed as integral
combinations of generators of the form given in (\ref{fundwedge2})
where
$$(|U|,|V|) = (1,2), (2,3) \mbox{ or } (3,4).$$

 As in the proof of Lemma~\ref{gener0}, we will show that each of these generators
$\rho_{U,V} \land
\gamma$ can be written as an integral  combination of generators given in (\ref{alpha2}),
which will complete the proof.

If $(|U|,|V|) = (1,2)$  then we are done.
If $(|U|,|V|) = (2,3)$ then we apply Lemma~\ref{gen0} since $m-2 +n-3 \equiv 0 \bmod 3$.
 Since  $m < n$, we have $m-2 \le n-3 \le 2(m-2) -3$.
Hence by Lemma~\ref{gen0}, we have that
$\rh_{\nu_{m,n}- |U|}(M_{[m]\setminus U,[n]\setminus V})$ is generated by wedge products
containing $\alpha$-cycles.  It follows that $\gamma$, and hence
$\rho_{U,V} \land
\gamma$,  is an integral  combination of wedge products containing $\alpha$-cycles.

Now suppose  $(|U|,|V|) = (3,4)$.
 Since $m < n$, we have $m-3 \le n-4 \le 2(m-3) -2$.   We can therefore apply
Lemma~\ref{gener1} since
$m-3 +n-4 \equiv 1
\bmod 3$. Hence,  $\rh_{\nu_{m,n}-
|U|}(M_{[m]\setminus U,[n]\setminus V})$ is generated by wedge products which contain
$\alpha$-cycles.
 It follows that $\gamma$, and hence
$\rho_{U,V} \land
\gamma$,  is an integral combination of wedge products containing $\alpha$-cycles.
\end{proof}

The next result follows readily from Lemma~\ref{gener2} by induction.

\begin{lem} \label{gen2} Suppose  $m+n \equiv 2 \bmod 3$ and $ m \le n \le 2m-4$. Then
$\rh_{\nu_{m,n}}(M_{m,n})$ is generated by elements of the form
$$\omega \land \gamma,$$ where $$\omega \in \rh_{\nu_{4,4}}(M_{U,V}), \qquad
\gamma \in \rh_{\nu_{m-4,n-4}}(M_{[m]\setminus U,[n]
\setminus V}),$$ and
$$ 4=|U| = |V|. $$
\end{lem}

\begin{thm} \label{chess2} Suppose $m+n \equiv 2 \bmod 3$ and $ m \le n \le 2m-13
$.  Then $\rh_{\nu_{m,n}}(M_{m,n})$ is a nontrivial
$3$-group of exponent at most 9.
\end{thm}

\begin{proof} Since $m-4 + n-4 \equiv 0 \bmod 3$ and $m-4 \le n-4 \le 2(m-4) - 9$, the result
follows from Lemma~\ref{gen2} and Theorem~\ref{chess0}.
\end{proof}

This completes the proof  of Theorem~\ref{ctor}.
We conjecture that the exponent  in Theorem~\ref{ctor} is $3$.  The following result shows
that  this conjecture need only be verified for
$m=n=9$.

\begin{thm} \label{conj99} For all $m,n$ that satisfy the hypothesis of
Theorem~\ref{ctor}, the exponent of $\rh_{\nu_{m,n}}(M_{m,n})$ divides the exponent of
$\rh_{\nu_{9,9}}(M_{9,9})$. Consequently if $\rh_{\nu_{9,9}}(M_{9,9})$ is an elementary
$3$-group then so is
$\rh_{\nu_{m,n}}(M_{m,n})$ for all $m,n$ that satisfy the hypothesis of
Theorem~\ref{ctor}.
\end{thm}

\begin{proof} The proof is similar to that of Theorem~\ref{chess2}.  It follows from
Lemmas~\ref{gen2} and \ref{gener0}.
\end{proof}

\subsection{Finite homology}

This subsection contains some partial results on the finite $\rh_{\nu_{m,n}}(M_{m,n})$
not covered by  Theorem~\ref{ctor}.
We start with an  analog of Corollary~\ref{3-prime-m}.

\begin{thm} \label{3-prime-c} The  Sylow 3-subgroup of
$\rh_{\nu_{m,n}}(M_{m,n})$  is nontrivial for all $m,n$ such that
$\rh_{\nu_{m,n}}(M_{m,n})$ is finite.
\end{thm}

\begin{proof} The proof is similar to that of (\ref{conjeq2}).  Assume $m \le n$ and
$\rh_{\nu_{m,n}}(M_{m,n})$ is finite with exponent $e$.

{\bf Case 1.} $m+n \equiv 1 \bmod 3$.  It follows from Theorem~\ref{FH}, that this case is
covered by Theorem~\ref{ctor} (i).

{\bf Case 2.} $m+n \equiv 0 \bmod 3$.   Consider the cycle
$z $  in the proof of (\ref{conjeq2}).  Recall  that $z$ cannot be a boundary in
$M_{[m] \uplus[n]^\prime}$. Since
$e z$ is a boundary in $M_{m,n}$, it is  also a boundary in $M_{[m] \uplus
[n]^\prime}$.  Since by Theorem~\ref{FH},  $7 \le m
\le n
\le 2m-6$, we have that
$m+n
\ge 15$. Therefore   Theorem~\ref{mtor} implies that $3 $ divides $e$, which means that
$\rh_{\nu_{m,n}}(M_{m,n})$ has 3-torsion.

{\bf Case 3.} $m+n \equiv 2 \bmod 3$.  By Theorem~\ref{FH}, we have $9 \le m \le n
\le 2m-7$. Consider the surjection
$\psi$ of Lemma~\ref{chessonto} (ii).  Since $m+n-4 \equiv 1 \bmod 3$ and $5 \le m-2 \le n-2
\le 2(m-2) - 5$, the range of $\psi$  has $3$-torsion by
Theorem~\ref{ctor} (i).  Since the domain is finite, the domain must also have 3-torsion.
\end{proof}

We have  not yet been able to eliminate $p$-torsion in finite $\rh_{\nu_{m,n}}(M_{m,n})$ for
primes
$p\ne 3$ except in the cases covered by Theorem~\ref{ctor}.  However, the lemmas of the
previous subsections provide an approach to  doing so  as well as   to reducing the exponent
in  Theorem~\ref{ctor} to $3$.  This approach, which  depends only on  anticipated improvements
of the computer computations, is  demonstrated by the following  result.

\begin{thm} \label{genmain}
\hspace{2in}
\begin{enumerate}
\item[(i)] If $m+n \equiv  0 \bmod 3$ and   $7 \le m \le n \le 2m-6$ then
$\rh_{\nu_{m,n}}(M_{m,n})$ is finite and its exponent divides the exponent of
$\rh_{\nu_{7,8}}(M_{7,8})$.
\item[(ii)] If $m+n \equiv 2 \bmod 3$ and  $11 \le m \le n \le 2m-10$ then
$\rh_{\nu_{m,n}}(M_{m,n})$ is finite and its exponent divides the exponent of
$\rh_{\nu_{7,8}}(M_{7,8})$.
\item[(iii)] If $m+n \equiv  2 \bmod 3$ and  $9 \le m \le n \le 2m-7$ and $ (m,n) \ne (10,10)
$ then $\rh_{\nu_{m,n}}(M_{m,n})$ is finite and  its exponent divides the
 exponent of $\rh_{\nu_{9,11}}(M_{9,11})$.
\end{enumerate}  Consequently if the Sylow $3$-subgroup of $\rh_{\nu_{7,8}}(M_{7,8})$ is
elementary then $\rh_{\nu_{m,n}}(M_{m,n})$ is an  elementary $3$-group for all
$m,n$ that satisfies the hypothesis of Theorem~\ref{ctor}.
\end{thm}

\begin{proof} Finiteness of the homology groups follow from Theorem~\ref{FH}.

(i)      We prove this by induction on $m$. The base case, $(m,n)=(7,8)$,
is trivial.  Now assume $m > 7$.  By Lemma~\ref{gener0}, the exponent
of $\rh_{\nu_{m,n}}(M_{m,n})$ divides the exponent of
$\rh_{\nu_{m-1,n-2}}(M_{m-1,n-2})$ if  $\rh_{\nu_{m-1,n-2}}(M_{m-1,n-2})$ is finite.   If
$m < n$ then $7 \le m-1 \le n-2 \le 2(m-1)-6$.  Hence by induction,
$\rh_{\nu_{m-1,n-2}}(M_{m-1,n-2})$ is finite and  the exponent of
$\rh_{\nu_{7,8}}(M_{7,8})$ is divisible by the exponent of
$\rh_{\nu_{m-1,n-2}}(M_{m-1,n-2})$ which is  divisible by the exponent of
$\rh_{\nu_{m,n}}(M_{m,n})$. If $m = n$ then $7 \le n-2 \le m-1 \le 2(n-2) - 6$.  So we can
apply the induction hypothesis in this case as well.

(ii)  By Lemma~\ref{gen2},  $\rh_{\nu_{m,n}}(M_{m,n})$ divides the exponent of
\newline
$\rh_{\nu_{m-4,n-4}}(M_{m-4,n-4})$ if $\rh_{\nu_{m-4,n-4}}(M_{m-4,n-4})$ is finite.
Since
$7
\le  m-4 \le n-4 \le 2(m-4) -6$, we can apply (i).

(iii) This is similar to the proof of (i) and is left to the reader.
\end{proof}

\begin{rem}{\rm We conjecture that there is some $m_0 $, such that if $n_0 = 2m_0 -6 $ or $n_0 =
2m_0-7 $ then
$\rh_{\nu_{m_0,n_0 }}(M_{m_0,n_0})$ is an elementary
$3$-group.  If this is so, then an argument like the one used in the proof of
Theorem~\ref{genmain} would yield the conclusion that
$\rh_{\nu_{m,n}}(M_{m,n})$ is an elementary $3$-group for all but a finite number of
pairs $(m,n)$ satisfying $m \le n \le 2m-5$.  (Recall $\rh_{\nu_{m,n}}(M_{m,n})$ is infinite
when
$n > 2m-5$..)}
\end{rem}

\section{Top homology of the chessboard complex} \label{top}

In this section we  construct bases for the top homology and cohomology of the
chessboard complex.   The basis for homology yields the
 decomposition result  used  in proving the torsion results of
Section~\ref{ctorsec}.

We assume familiarity with  the  representation theory of the symmetric group $\Ss_n$ and
tableaux combinatorics, cf.,
\cite{Sa},
\cite{St},
\cite{F}.  The Specht module   (or irreducible representation   of $\Ss_n$ ) over
$\C$
 indexed by the partition  $\lambda \vdash n$, is denoted by $S^\lambda$.   Recall  that
the dimension of $S^\lambda$ is the number $f^\lambda$ of standard Young tableaux of shape
$\lambda$.

The direct product $\Ss_m \times \Ss_n$ acts  on the chessboard
complex $M_{m,n}$ by relabelling the graph vertices in $[m]$ and
$[n]^\prime$, and this induces a representation of $\Ss_m \times
\Ss_n$ on $\tilde H_*(M_{m,n};\C)$.  The following result enables
one to express the Betti numbers in terms of the number of pairs
of standard Young tableaux of certain shapes.

\begin{thm} [Friedman and Hanlon \cite{FH}] \label{FH2}  For all $p,m,n \in \Z$,
where $m,n \ge 1$, the following isomorphism of
$(\Ss_m
\times
\Ss_n)$-modules holds:
$$\tilde H_{p-1}(M_{m,n};\C)  \cong_{\Ss_m \times \Ss_n} \bigoplus_{(\lambda,\mu)
\in \RR(m,n,p)}
 S^{\lambda^\prime} \otimes S^\mu,$$ where
$\RR(m,n,p)$ is the set of all pairs of partitions $(\lambda \vdash m, \mu \vdash n)$
that can be obtained in the following
way.  Take a partition $\nu \vdash p$  that contains an $(m-p) \times (n-p)$ rectangle
but contains no $(m-p+1) \times
\times ( n-p+1)$ rectangle.  Add a column of size $m-p$ to $\nu$ to obtain $\lambda$ and
add a
row of size $n-p$ to $\nu$ to obtain $\mu$.
 See Figure~\ref{top}.1.

\begin{center}
\includegraphics[width=3.5in]{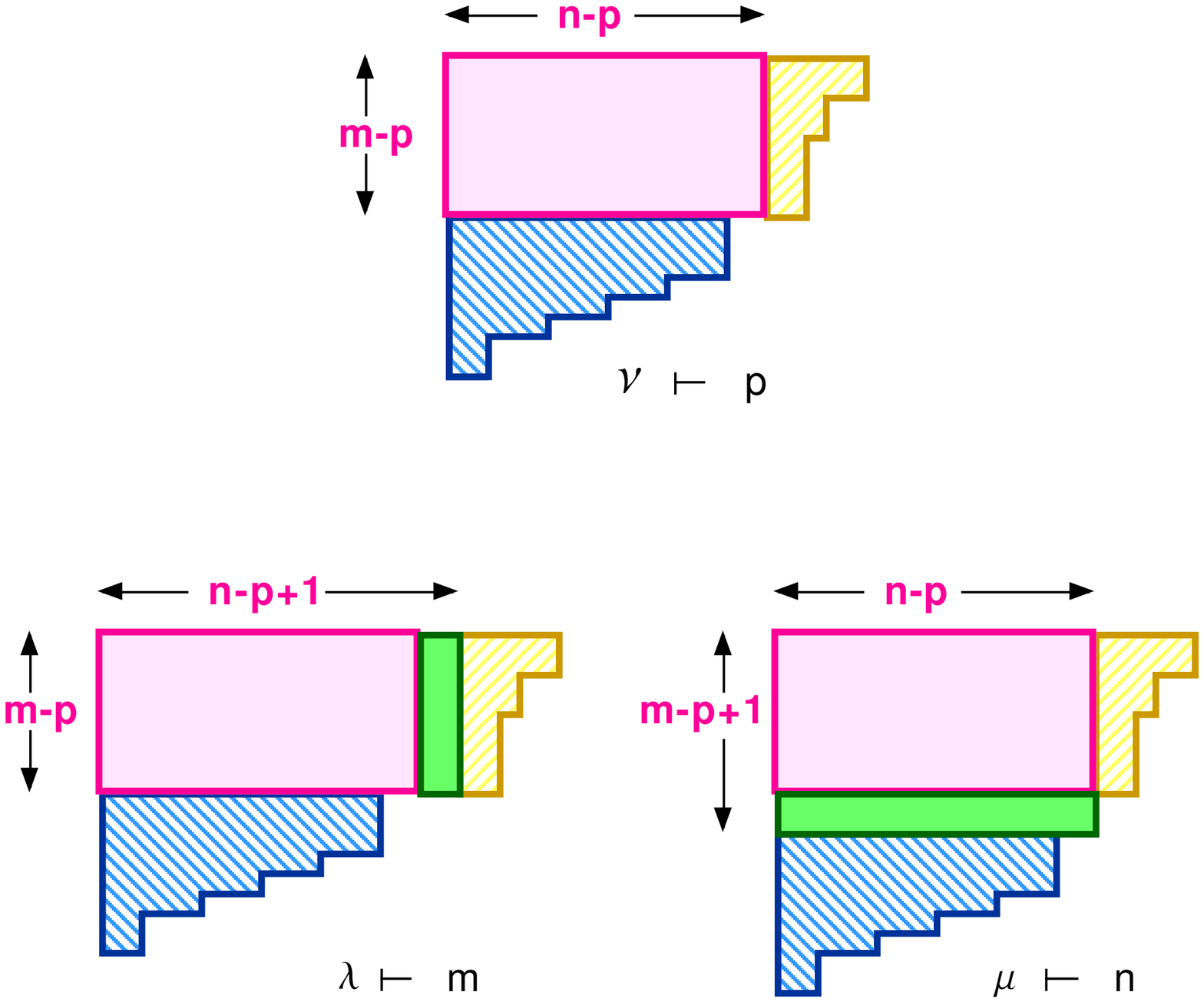}

{\bf Figure~\ref{top}.1}  \end{center}
\end{thm}

\begin{cor}[Garst\cite{G}] \label{garstcor} For all $m \le n$, the following isomorphism of
$\mathfrak S_n$-modules holds
$$\tilde H_{m-1}(M_{m,n};\C)  \cong_{ \Ss_n} \bigoplus_{\scriptsize\begin{array}{c}
\lambda
\vdash m
\\
\lambda_1 \le n-m \end{array} }
 f^{\lambda} \, S^{\lambda^*},$$
where   $\lambda^*$
is the partition obtained from
$\lambda$ by adding a part of size $n-m$.
\end{cor}

It follows immediately from  Corollary~\ref{garstcor} that the rank of the top homology
$\rh_{m-1}(M_{m,n})$ of the chessboard complex
$M_{m,n}$ is the number of pairs  of standard Young tableaux $(S,T)$ such that $S$ has $m$
cells, $T$ has $n$ cells and the shape of $S$ is the same as the shape of
$T$ minus the first row.   Let ${\mathcal P}_{m,n}$ be the set of such pairs of standard
tableaux.   We construct for each $(S,T)
\in {\mathcal P}_{m,n}$, a cycle $\eta(S,T) \in \rh_{m-1}(M_{m,n})$, and show that these
cycles form a basis for homology.

In order to prove that the $\eta(S,T)$ form a
basis for homology, we  construct  cocycles $\gamma(S,T) $ which form a basis for
cohomology.  Since our complex is finitely generated we can view the cohomology group as a
subquotient of the chain group, just as is done for the homology group.  Indeed, for any finite
simplicial complex
$\Delta$ on  vertex set
$\{x_1,\dots,x_r\}$, let
$\langle
\phantom{a},\phantom{b}
\rangle$ be the bilinear form on
$C_{k-1}(\Delta)$ for which the oriented simplices $(x_{i_1},\dots,x_{i_k})$,
 $i_1 <\dots < i_k$, form an orthonormal basis.  The coboundary map $\delta_k:C_k(\Delta) \to
C_{k+1}(\Delta)$ is the adjoint of the boundary map.  That is
$$\langle u,\delta_k(v) \rangle = \langle \partial_{k+1}(u),v\rangle , $$
for all $u \in C_{k+1}(\Delta)$ and $v \in C_k(\Delta)$. The $k${\em th} cohomology group is
defined to be the quotient of the cocycle group $Z^k(\Delta):= \ker \delta_k$ by the coboundary
group $B^{k}(\Delta):= \im \delta_{k-1}$.

Our construction of the cycles and cocycles uses the classical
Robinson-Schensted correspondence.  We begin with the cocycles.
Let $(S,T) \in \mathcal P_{m,n}$.  First add a cell with entry
$\infty$ to the bottom of each  of the first $n-m$ columns (some
may be empty) of $S$
 to obtain a semistandard tableau $S^*$ of the same shape as $T$.  (Here $\infty$ represents
a number larger than $m$.) See Figure~\ref{top}.2.  The inverse of the Robinson-Schensted
bijection applied to
$(S^*,T)$ produces a permutation $\sigma$ of the multiset
$\{1,2,\dots,m,\infty^{n-m}\}$.  The multiset permutation $\sigma$ corresponds naturally to
the oriented simplex  of
$M_{m,n}$ given by
\beq \label{word}\tau(\sigma):= \left(\sigma(i_1)i_1^\prime, \sigma(i_2)i_2^\prime,\dots,
\sigma(i_m)i_m^\prime\right),\eeq
 where
$\sigma(i_1)\sigma(i_2)\cdots\sigma(i_m)$ is the subword of $\sigma=
\sigma(1)\sigma(2) \cdots \sigma(n)$ obtained by removing the
$\infty$'s. This oriented simplex is clearly a cocycle since it is in the top dimension.  Let
$\gamma(S,T)$ be the coset of  the coboundary group
$B^{m-1}(M_{m,n})$ that contains this oriented simplex.

We demonstrate the procedure for constructing $\gamma(S,T)$ by letting $(S,T)$ be the pair of
tableaux given in Figure~\ref{top}.2. After applying the inverse of Robinson-Schensted to
$(S^*,T)$ we have the multiset permutation
$$\infty\,\infty \,2\,\infty \,4
\,\infty\, 3
\,1.$$  The oriented simplex that corresponds to this multiset permutation is
$$(23^\prime,45^\prime,37^\prime,18^\prime).$$  Hence,
$\gamma(S,T)$ is  the coset of $B^{3}(M_{4,8})$ that contains the oriented simplex
$(23^\prime,45^\prime,37^\prime,18^\prime)$.

\vspace{.2in}
\begin{center}
\includegraphics[width=2.2in]{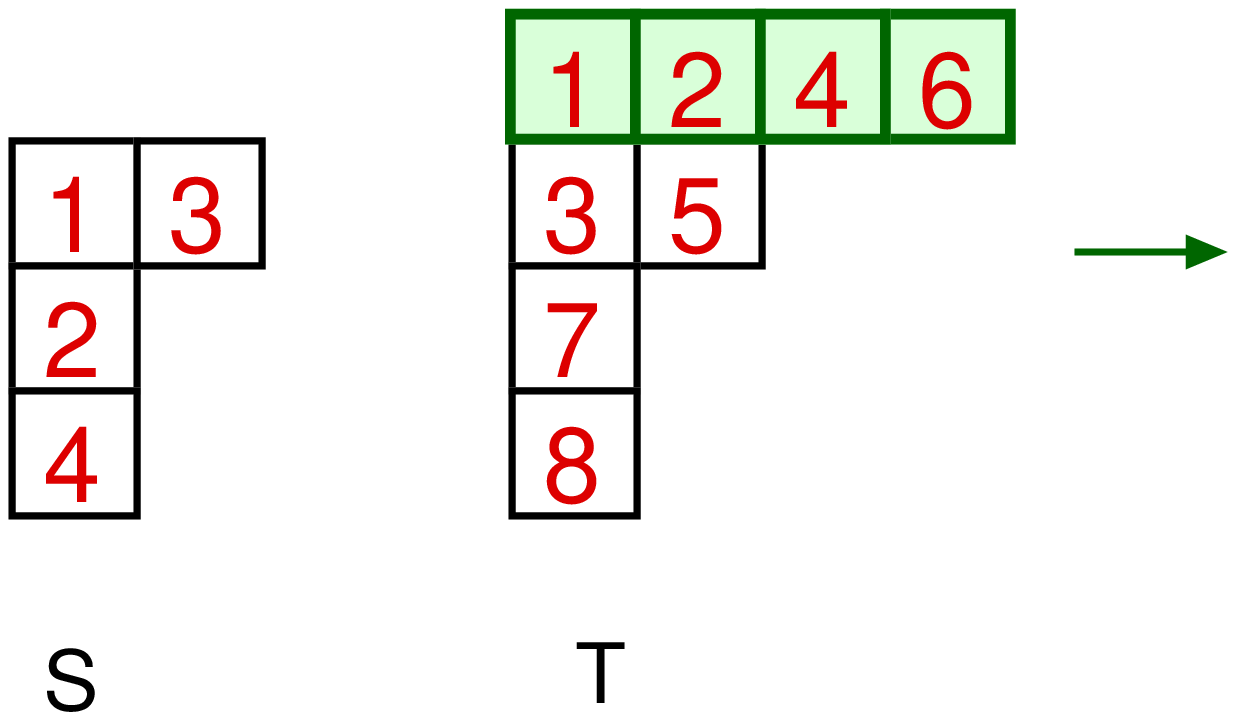}
\includegraphics[width=2.2in]{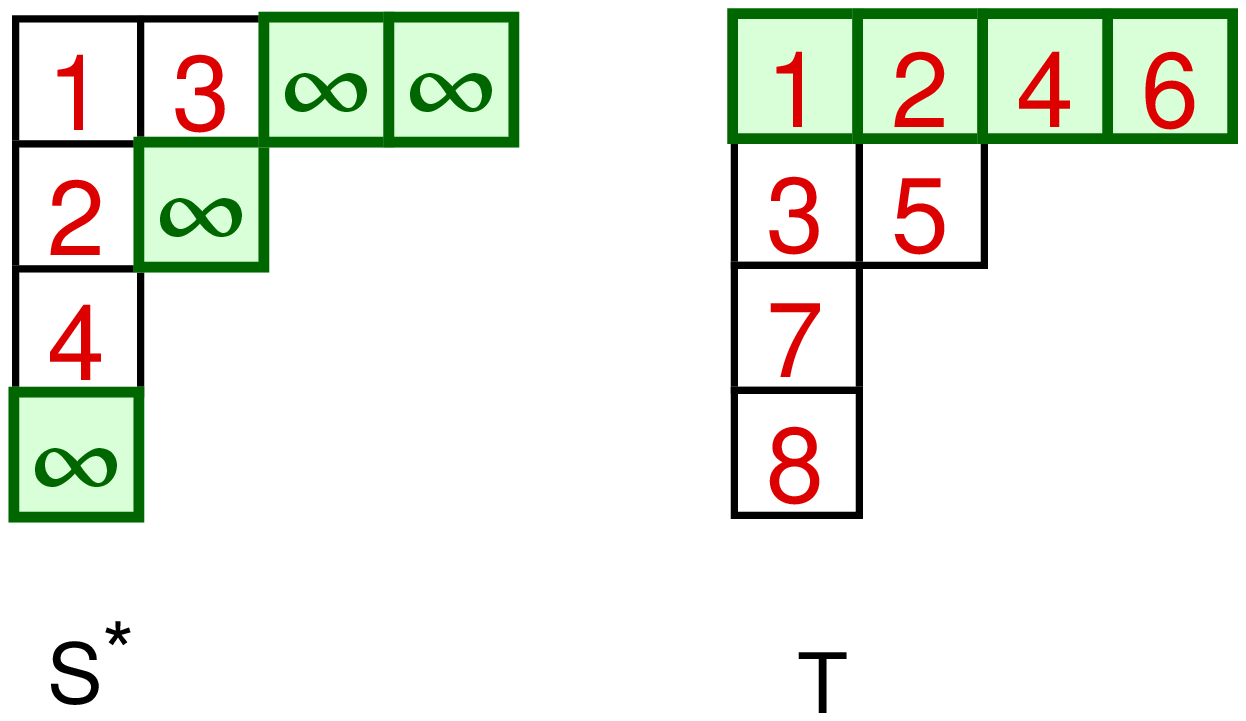}
\vspace{.2in} {\bf Figure \ref{top}.2 }\end{center}

The construction of the cycles is a bit more involved.  Recall that in  the inverse
Robinson-Schensted procedure,  an entry ``pops'' from a  cell in the top row of the left
tableau when an entry is ``crossed out'' of the right tableau.  For each top cell, we must
keep track  of the entries of $S^*$ that are popped and the corresponding entries of $T$ that
are crossed out.  For each $i = 1,2,\dots, n-m$, let $A_i^*$ be the multiset of entries that
are popped from the
$i$th
 cell
 of the top row of
$S^*$ and let
$B_i$ be the corresponding set of entries that are crossed out of $T$.  One can easily see
that $A^*_i$ is actually a set and
$\infty \in A^*_i$ for all $i$.  Now let $A_i = A^*_i\setminus \{\infty\}$.  So
$|A_i| = |B_i|-1$.  It is easily observed that
$M_{A,B}$ is an orientable pseudomanifold whenever $|A| =|B|-1$, which implies that its top
homology  is cyclic. The fundamental cycle of $M_{A,B}$ (that is, generator of top homology,
which is unique up to sign) is explicitly given by
\beq \label{fund} \rho_{A,B} := \sum_{\sigma \in\mathfrak
S_{A\cup\{\infty\}}}
\sgn(\sigma) \tau(\sigma).\eeq
  Now define
$$\eta(S,T) = \rho_{A_i,B_i} \land \cdots \land \rho_{A_{n-m},B_{n-m}}.$$

We demonstrate the procedure for constructing $\eta(S,T)$ on the tableaux $S,T$ of
Figure~\ref{top}.2.  Refer to Figure~\ref{top}.3.  First  entry $8$  is crossed out of $T$ and
entry
$1$ is popped from the first cell of the first row of $S^*$.  So $1$ is placed in  $A^*_1$  and
$8$ is placed in  $B_1$.  Next entry $7$ is crossed out and entry $3$ is popped from the
second cell.  So $3$ is placed in $A^*_2$ and
$7$ is placed in $B_2$.  We eventually end up with
$$A^*_1=\{1,2,\infty \},\,  A^*_2 =\{3,4, \infty \}, \, A^*_3 = A^*_4 = \{\infty\},$$
$$B_1 = \{1,3,8\},\, B_2 = \{2,5,7\}, \,B_3 = \{4\}, \,B_4 = \{6\}.$$
 Hence
 $$A_1=\{1,2\},\,  A_2 =\{3,4\}, \, A_3 = A_4 = \emptyset$$ Now $$\eta(S,T) =
\rho_{\{1,2\},\{1,3,8\}} \land \rho_{\{3,4\},\{2,5,7\}}.$$

\vspace{.02in}
\begin{center}
\includegraphics[width=2.7in]{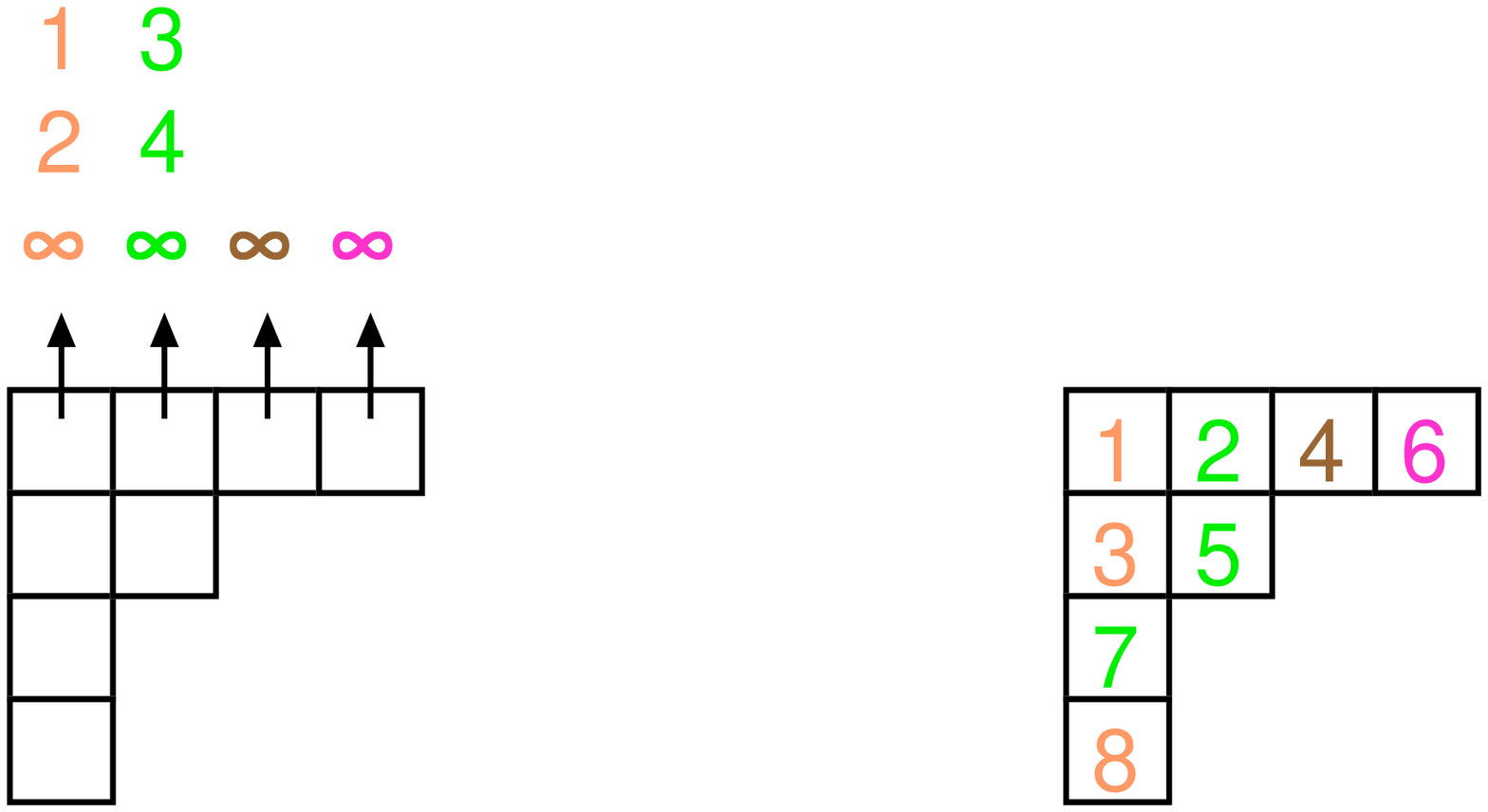}

\vspace{.1in}
{\bf Figure \ref{top}.3 }\end{center}

\vspace{.2in}

\begin{thm} \label{basis} Let $m \le n$.  Then
\begin{itemize}
\item
$\{\eta(S,T) : (S,T) \in  {\mathcal P_{m,n}} \}$ is a basis for
$\rh_{m-1}(M_{m,n})$.
\item
$\{\gamma(S,T) : (S,T) \in  {\mathcal P_{m,n}} \}$ is a basis for the free part of
$\rh^{m-1}(M_{m,n})$.
\end{itemize}
\end{thm}

We  need some general theory in order to prove this result.  For
any abelian group $G$, let $ G_{\mbox{tor }}$ denote the subgroup
of $G$ consisting of torsion elements of $G$

\begin{prop} \label{pairing} Let $\Delta$ be a simplicial complex. Suppose
\begin{itemize}
\item $r = \mbox{rank}(\tilde H_k(\Delta)/ \tilde H_k(\Delta)_{\mbox{tor }})$,
\item
$u_1,\dots,u_r \in
Z_k(\Delta)$,
\item $v_1,\dots,v_r \in  Z^k(\Delta)$,
\item the matrix $(\langle
u_i,v_j \rangle)_{i,j = 1 \dots,r}$ is invertible over $\Z$.
\end{itemize}
Then $\{\hat u_1,
\dots,
\hat u_r \}$ is a basis for $\tilde H_k(\Delta)/\tilde H_k(\Delta)_{\mbox{tor }} $ and $\{\hat
v_1, \dots,
\hat v_r \}$ is a basis for $\tilde H^k(\Delta)/\tilde H^k(\Delta)_{\mbox{tor } }$,
where $\hat x$ denotes the coset of $\tilde H^k(\Delta)_{\mbox{tor } }$ or $ \tilde
H_k(\Delta)_{\mbox{tor }}$ containing $\bar x$.
\end{prop}

\begin{proof} The invertibility of the matrix $A:= (\langle
u_i,v_j \rangle)_{i,j = 1 \dots,r}$ implies that $\bar u_1,\dots,\bar u_r$ are
independent in $\rh_{k}(\Delta,\Q)$.
Since $r = \dim
\rh_{k}(\Delta,\Q)$, we have that $\bar u_1,\dots,\bar u_r$ also spans $\rh_{k}(\Delta,\Q)$.

Let $u \in Z_k(\Delta)$.  Then $u \in Z_k(\Delta,\Q)$.  So
$$\bar u = \sum_{i=1}^r c_i \, \bar u_i, \quad c_i \in \Q$$
in $\rh_k(\Delta,\Q)$.   This means
$$u - \sum_{i=1}^r c_i \,  u_i = \partial(y)$$
for some $y \in C_{k+1}(\Delta,\Q)$.
For each $j$,  we have
$$\langle u, v_j \rangle - \sum_{i=1}^r c_i\langle u_i,v_j\rangle = \langle \partial(y),
v_j\rangle = \langle y, \delta(v_j) \rangle = 0,$$  since $v_j$ is
a cocycle. It follows that
$$\left[\begin{array}{c}\langle u,v_1\rangle \\ \vdots \\ \langle u,v_r\rangle
\end{array}\right] = A \left[\begin{array}{c}c_1 \\ \vdots \\c_r
\end{array}\right] ,$$
which implies
$$\left[\begin{array}{c}c_1 \\ \vdots \\c_r
\end{array}\right] = A^{-1} \left[\begin{array}{c}\langle u,v_1\rangle \\ \vdots \\ \langle
u,v_r\rangle
\end{array}\right] \in \Z^r.$$

Let $t \in \Z^+$ be such that $ty \in C_{k+1}(\Delta)$.  Since $$t(u - \sum_{i=1}^r c_i u_i) =
\partial(ty),$$  we have $\bar u - \sum_{i=1}^r c_i \bar u_i  \in  \rh_k(\Delta)_{\mbox{tor}}$.
  It
follows that $$\hat u = \sum_{i=1}^r c_i \hat u_i$$ in $\rh_k(\Delta) /
\rh_k(\Delta)_{\mbox{tor}}$.  Hence $\hat u_1, \dots,\hat u_r$ generates $\rh_k(\Delta) /
\rh_k(\Delta)_{\mbox{tor}}$.  Since $r = \mbox{rank}(\rh_k(\Delta) /
\rh_k(\Delta)_{\mbox{tor}})$, these elements form a basis for $\rh_k(\Delta) /
\rh_k(\Delta)_{\mbox{tor}}$.  By symmetry $\hat v_1, \dots,\hat v_r$ forms a basis for
$\rh^k(\Delta) /
\rh^k(\Delta)_{\mbox{tor}}$.
\end{proof}

\begin{proof}[Proof of Theorem~\ref{basis}]
  For  $(S,T) \in  \mathcal P_{m,n}$, let
\beq \nonumber v(S,T) := \tau(\mbox{RS}^{-1}(S^*,T)) \in C_{m-1}(M_{m,n}),\eeq
where $\mbox{RS}^{-1}$ denotes the inverse of the Robinson-Schensted map and $\tau$ is the map
defined in (\ref{word}).  Let
\beq \nonumber u(S,T): = \sum_{\omega \in \mathfrak S_{B_1} \times \cdots \times
\mathfrak S_{B_{n-m}}}\,
\sgn(\omega)\,\,\tau(\mbox{RS}^{-1}(S^*,T)\omega) \in C_{m-1}(M_{m,n}),\eeq
where $B_1,\dots, B_{n-m} $ are the sets defined  in the construction of $\rho(S,T)$.

For all $(S,T) \in  {\mathcal P_{m,n}}$, we have
\beq \label{cobasis} \gamma(S,T) = \ov{v(S,T)},\eeq
where  $\ov {x}$
denotes the cohomology class of $x$ in $\tilde H^{m-1}(M_{m,n})$. It is not hard to
see that
\beq \label{hombasis} \eta(S,T)  = \sgn(B_1,\dots,B_{n-m})\,\,\, u(S,T),\eeq
where $\sgn(B_1,\dots,B_{n-m})$ is the sign of the permutation
obtained by concatenating the words obtained by writing  each $B_i$ in decreasing
order.

Next we claim that for all $(S_1,T_1), (S_2,T_2) \in \mathcal P_{m,n}$,
\beq \label{lex} \\ \nonumber \langle
u(S_1,T_1),
v(S_2,T_2)\rangle
\ne 0\quad  \Longrightarrow \quad  \mbox{RS}^{-1}(S_2^*,T_2) \le_{\mbox{lex}}
\mbox{RS}^{-1}(S_1^*,T_1)\eeq  where $\le_{\mbox{lex}}$ denotes  lexicographical order.
Note that the subword of $\mbox{RS}^{-1}(S^*,T)$ obtained by
restricting to the positions in $B_i$, is decreasing for each $i = 1, \dots, n-m$. Hence any
rearrangement of letters of $\mbox{RS}^{-1}(S^*,T)$ occupying positions in $B_i$, produces a
lexicographically smaller word.  Hence for each $\omega \in
\mathfrak S_{B_1}
\times
\cdots
\times
\mathfrak S_{B_{n-m}}-\{e\}$,
$$\mbox{RS}^{-1}(S^*,T)\omega \,\,<_{\mbox{lex}}\,\, \mbox{RS}^{-1}(S^*,T).$$
The claim (\ref{lex}) follows from this.
We also have  that
\beq \label{diag} \langle
u(S,T),
v(S,T)\rangle = 1\eeq for all $(S,T) \in \mathcal P_{m,n}$.

Now order  the pairs of standard
tableaux
$$(S_1,T_1),\dots,(S_r,T_r)$$ in
$\mathcal P_{m,n}$
so that $\mbox{RS}^{-1}(S_i^*,T_i) <_{\mbox{lex}}\mbox{RS}^{-1}(S_j^*,T_j)$
if $i <j$.  It follows from (\ref{lex}) and (\ref{diag})  that  the matrix
$$(\langle u(S_i,T_i), v(S_j,T_j)\rangle)_{i,j = 1,\dots,r}$$ is unitriangular.
There is no torsion in the top homology, and by Corollary~\ref{garstcor}, $|\mathcal P_{m,n}| =
\mbox{rank}
\rh_{m-1}(M_{m,n})$. Hence the result follows from (\ref{cobasis}), (\ref{hombasis}) and
Proposition~\ref{pairing}.
\end{proof}

\begin{cor}\label{topcor} Let $m \le n$.  Then $\rh_{m-1}(M_{m,n})$ is generated by cycles of
the form
$$\rho_{A,B} \land \tau,$$ where
\begin{itemize}
\item $A \subseteq [m]$, $B \subseteq [n]$ and   $1 \le |A| = |B| -1$
\item $\rho_{A,B}$ is a fundamental cycle of the pseudomanifold
$M_{A,B}$
\item $\tau \in \rh_{m-1-|A|}(M_{[m]-A,[n]-B})$.
\end{itemize}
\end{cor}

\section{Infinite homology of the chessboard complex} \label{infinsec}

In this section we study torsion in infinite $\tilde
H_{\nu_{m,n}}(M_{m,n})$.   Recall from Theorem~\ref{FH} that for $m \le n$,  the homology group
$\tilde H_{\nu_{m,n}}(M_{m,n})$ is infinite if and only if $n \ge 2m-4$ or $(m,n) \in \{(6,6),
(7,7),(8,9)\}$.  From Table 1.2, we see that there is 3-torsion if  $(m,n) =
(6,6)$ or $ (7,7)$.   We expect that there is 3-torsion for $(m,n) = (8,9)$ as well,  but
 have not yet been able to verify this by computer.

\begin{conj} \label{conjfree} Let $m \le n$. Then $\tilde H_{\nu_{m,n}}(M_{m,n})$ is free if
and only if
$n \ge 2m-4$.
\end{conj}

The conjecture clearly holds in the case that $n \ge 2m-1$, since in this case $\nu_{m,n}=
m-1$, which means that  $\tilde H_{\nu_{m,n}}(M_{m,n})$ is top homology.   The conjecture  for
$n = 2m-2$ is proved in the following result.  The cases $n = 2m-3$ and $n = 2m-4$ are left
open.

\begin{thm} \label{free} If $n = 2m-2$ then
$$\rh_{\nu_{m,n}}(M_{m,n}) \cong \Z^{c_{m-1}},$$
where $c_m$ is the Catalan number
$ \frac 1 {m+1} \binom {2m}{ m}$.
\end{thm}

\begin{proof}  Theorem~\ref{FH2} applied to $\rh_{\nu_{m,2m-2}}(M_{m,2m-2};\C)$ yields a
particularly nice formula.  First note that
$\nu_{m,2m-2} = m-2$.   Next  observe that the set
$\mathcal R(m,2m-2,m-1)$   consists of a single pair of partitions; namely the pair
$((m),(m-1)^2)$.  Hence  Theorem~\ref{FH2} yields,
$$\rh_{\nu_{m,2m-2}}(M_{m,2m-2};\C) \cong_{\mathfrak S_n} S^{(m-1)^2}.$$
It follows that the degree $\nu_{m,2m-2}$ Betti number of $M_{m,2m-2}$ is $f^{(m-1)^2}$, the
number of standard Young tableaux of shape $(m-1)^2$. Hence
\beq \\ \nonumber \label{rank}  \mbox{rank } (\rh_{\nu_{m,2m-2}}(M_{m,2m-2})/
\rh_{\nu_{m,2m-2}}(M_{m,2m-2})_{\mbox{tor }})
 = f^{(m-1)^2}.\eeq
Since the number of standard Young
tableaux of shape
$(m-1)^2$ is the Catalan number
$c_{m-1}$, we need only show that $\rh_{\nu_{m,2m-2}}(M_{m,2m-2})$ is free.

Given a partition $\lambda$, let $S^{\lambda}_{\Z}$ denote the Specht module indexed by
$\lambda$ with integer coefficients.  It is well-known that
$S^{\lambda}_{\Z}$ is a free group of rank $f^\lambda$, which is isomorphic to the group
generated by the
$\lambda$-tableaux subject to the column relations
 and the Garnir relations.  For $\lambda = (m-1)^2$,
these relations can be described as follows:
$$ \begin{array}{|c|c|c|}
\hline \cdots & a_j & \cdots \\
\hline \cdots & b_j & \cdots \\ \hline
\end{array}\,\,+\,\,  \begin{array}{|c|c|c|}
\hline \cdots & b_j & \cdots \\
\hline \cdots & a_j & \cdots \\ \hline
\end{array} $$

$$ \begin{array}{|c|c|c|c|}
\hline \cdots & a_{j-1} & a_{j}& \cdots \\
\hline \cdots & b_{j-1}&\cdot & \cdots \\ \hline
\end{array}\,\,- \,\,  \begin{array}{|c|c|c|c|}
\hline \cdots & a_{j} & a_{j-1}& \cdots \\
\hline \cdots & b_{j-1} &\cdot & \cdots \\ \hline
\end{array} \,\, + \,\,
 \begin{array}{|c|c|c|c|}
\hline \cdots & a_{j} & b_{j-1}& \cdots \\
\hline \cdots & a_{j-1}&\cdot & \cdots \\ \hline
\end{array} $$

$$ \begin{array}{|c|c|c|c|}
\hline \cdots &\cdot & a_{j}& \cdots \\
\hline \cdots & b_{j-1}&\ b_{j} & \cdots \\ \hline
\end{array}\,\,- \,\,  \begin{array}{|c|c|c|c|}
\hline \cdots & \cdot & a_{j}& \cdots \\
\hline \cdots & b_{j}& b_{j-1} & \cdots \\ \hline
\end{array} \,\, + \,\,
 \begin{array}{|c|c|c|c|}
\hline \cdots & \cdot & b_{j}& \cdots \\
\hline \cdots & a_{j}&b_{j-1} & \cdots \\ \hline
\end{array}$$

Let
$\phi: S^{(m-1)^2}_\Z \to \rh_{\nu_{m,2m-2}}(M_{m,2m-2})$
be the homomorphism defined on generators by

$$\phi\left ( \,\,\,\begin{array}{|c|c|c|c|}
\hline a_1 & a_2 & \cdots & a_{m-1} \\
\hline b_1 & b_2 & \cdots & b_{m-1} \\ \hline
\end{array}\,\,\, \right ) = \alpha_{1,a_1^\prime,b_1^\prime} \land
\alpha_{2,a_2^\prime,b_2^\prime} \land \dots \land
\alpha_{m-1,a_{m-1}^\prime,b_{m-1}^\prime}.$$
 To verify that this map is well defined we need only check that the three relations
for the Specht module given above are mapped to $0$ in $\rh_{\nu_{m,2m-2}}(M_{m,2m-2})$.
For the first relation we have $$ \phi\left ( \,\,\,\begin{array}{|c|c|c|}
\hline \cdots & a_j & \cdots \\
\hline \cdots & b_j & \cdots \\ \hline
\end{array}\,\,+\,\,  \begin{array}{|c|c|c|}
\hline \cdots & b_j & \cdots \\
\hline \cdots & a_j & \cdots \\ \hline
\end{array}\,\,\,\right) \,\, = \,\, \dots \land (\alpha_{j,a^\prime_j,b^\prime_j} +
\alpha_{j,b^\prime_j,a^\prime_j})\land
\dots,$$
which is clearly $0$.

For the second relation, we have $$ \phi\left ( \,\,\, \begin{array}{|c|c|c|c|}
\hline \cdots & a_{j-1} & a_{j}& \cdots \\
\hline \cdots & b_{j-1}&\cdot & \cdots \\ \hline
\end{array}\,\,- \,\,  \begin{array}{|c|c|c|c|}
\hline \cdots & a_{j} & a_{j-1}& \cdots \\
\hline \cdots & b_{j-1} &\cdot & \cdots \\ \hline
\end{array} \,\, + \,\,
 \begin{array}{|c|c|c|c|}
\hline \cdots & a_{j} & b_{j-1}& \cdots \\
\hline \cdots & a_{j-1}&\cdot & \cdots \\ \hline
\end{array}\,\,\,\right) $$

\begin{eqnarray*} \quad =& \dots & \land \,\,\,
\big((\alpha_{j-1,a^\prime_{j-1},b^\prime_{j-1}}
\land
\alpha_{j,a^\prime_{j},b^\prime_{j}})\, -\,(\alpha_{j-1,a^\prime_{j},b^\prime_{j-1}}
\land
\alpha_{j,a^\prime_{j-1},b^\prime_{j}})
\\  &+&\!\!\!(\alpha_{j-1,a^\prime_{j},a^\prime_{j-1}} \land
\alpha_{j,b_{j-1}^\prime,b^\prime_{j}})
\big)\,\,\,\land
\quad \dots .
\end{eqnarray*}
We will  show that this cycle, which we denote by $\rho$, is  a boundary.
After cancelling terms we get \begin{eqnarray*}\rho = & \dots & \land \,\,\,
\big((\alpha_{j-1,a^\prime_{j-1},b^\prime_{j-1}}
\land ja_{j}^\prime)\, -\,(\alpha_{j-1,a^\prime_{j}, b^\prime_{j-1}} \land
ja_{j-1}^\prime)
\\&+&\!\!\!(\alpha_{j-1,a^\prime_{j},a^\prime_{j-1}} \land
jb_{j-1}^\prime)\big)\,\,\,\land
\quad \dots,\end{eqnarray*}
which is   an element of the chain group $C_{m-2}(M_{[m-1],[2m-2]\setminus \{b_j\}})$.
Hence $mb_{j}^\prime\land\rho \in C_{m-1}(M_{m,2m-2})$.
Since
$\partial(mb_{j}^\prime
\land
\rho) = \rho$, the second relation maps to 0.
By symmetry the third relation maps to 0 as well. Hence $\phi$ is a well defined homomorphism.

We claim that $\phi$ is surjective.  Indeed, by Lemma~\ref{gen1},
$\rh_{\nu_{m,2m-2}}(M_{m,2m-2})$ is generated by elements of the form  $$
\alpha_{\sigma(1),a_1^\prime,b_1^\prime} \land
\alpha_{\sigma(2),a_2^\prime,b_2^\prime} \land \dots \land
\alpha_{\sigma(m-1),a_{m-1}^\prime,b_{m-1}^\prime}.$$
It follows from Lemma~\ref{rel} that $\sigma$ can be taken to be the identity permutation,
which means that $\rh_{\nu_{m,2m-2}}(M_{m,2m-2})$ is generated by the images of the
$(m-1)^2$-tableaux.

Let $$\pi: \rh_{\nu_{m,2m-2}}(M_{m,2m-2}) \to \rh_{\nu_{m,2m-2}}(M_{m,2m-2})/
\rh_{\nu_{m,2m-2}}(M_{m,2m-2})_{\mbox{tor
}},$$ be the projection map. The composition
$$\pi \circ \phi:S^{(m-1)^2}_\Z \to \rh_{\nu_{m,2m-2}}(M_{m,2m-2})/
\rh_{\nu_{m,2m-2}}(M_{m,2m-2})_{\mbox {tor
}},$$ is a surjective homomorphism between free
groups.  Since these groups have equal rank by (\ref{rank}),  the composition   $\pi \circ
\phi$ is an isomorphism, which implies that the surjection $\phi$ is an isomorphism as well.
We can now conclude that $\rh_{\nu_{m,2m-2}}(M_{m,2m-2})$ is free. \end{proof}

\begin{cor} The set $$\{\phi(T) : T \,\,{\rm a
\,\,standard\,\, tableau\,\, of\,\, shape}\,\, (m-1)^2\}
$$ is a basis for $\tilde H_{\nu_{m,2m-2}}(M_{m,2m-2})$.
\end{cor}

In \cite[Section 9.1]{BBLSW}, it is observed that when $m=n$, $M_{m,n}$ collapses to an
$(n-2)$-dimensional complex. Hence for $m=n$, the homology group $\rh_{i}(M_{m,n})$ is free
whenever  $i \ge m-2$.
  Theorem~\ref{free} implies that the same is true for $n=2m-2$.  This and the computer data
suggest the following  conjecture, which implies Conjecture~\ref{conjfree}.

\begin{conj} Let $m \le n$ and $i \ge \nu_{m,n}$.  Then $\rh_{i}(M_{m,n})$ is free if and only
if
$i
\ge m-2$.
\end{conj}

\section{Subcomplexes of the chessboard complex} \label{subsec}

Our  goal in this section is to  establish sharpness of  a connectivity bound for  the
simplicial complex of nontaking rooks on an
$n\times n$ chessboard with a diagonal removed.  This bound was obtained by Bj\"orner
and Welker
\cite{BW} as a consequence of a more general result of Ziegler \cite{Z} on nonrectangular
boards.

For any subset $A$ of the set of positions on an $m \times n$ chessboard, let $M(A)$ be the
simplicial complex of nontaking rooks on  $A$.  That is, for $A \subseteq [m] \times [n]$, the
simplicial complex
$M(A)$ has vertex set
$A$ and faces
 $\{(i_i, j_1), (i_2,j_2), \dots, (i_k,j_k)\}\subseteq A$ such that $i_s \ne i_t$ and $j_s
\ne i_t$ for all
$s
\ne t$.
Let $$D_n= [n] \times [n] \setminus
\{(1,1),(2,2), \dots, (n,n) \}.$$

\begin{thm} [Bj\"orner and Welker \cite{BW}]\label{diagth}  For all $n \ge 2
$, the simplicial complex
$M(D_n)$ is
$(\nu_{2n}-1)$-connected.
\end{thm}

 Bj\"orner and Welker
\cite{BW} use computer calculations to obtain the following table which establishes
sharpness of their connectivity bound for $3\le n\le7$.  We will use results of the previous
sections to establish sharpness for $n > 7$.

$$\begin{array}{|c||c|c|c|c|c|c|c|}
\hline  n    & 2 & 3&4&5&6&7\\ \hline
  \rh_{\nu_{2n}}(D_n) &  0  & \Z^2 & \Z^4 & \Z & \Z^{24} \oplus \Z^5_3 & \Z^{415} \oplus
\Z^{15}_3\\
\hline
\end{array}$$

\vspace{.1in}
\begin{center}
Table 8.1
\end{center}

For $n \ge 3$ and $i= 0, \dots, \lfloor {n \over 3} \rfloor -1$, let $$S_i =
\{(3i+1,3i+1),(3i+1, 3i+2), (3i+2, 3i+3),(3i+3,3i+3)\},$$
and let $$B_n = (\biguplus_{i=0}^N S_i) \,\, \uplus\,\, R_n,$$
where
$$N = \begin{cases}{n-3 \over 3}   &\mbox{ if  }n \equiv 0 \bmod 3 \\
{n-5 \over 3}  &\mbox{ if  }n \equiv 2 \bmod 3\\
{n-7 \over 3} &\mbox{ if  }n \equiv 1 \bmod 3\end{cases}$$
and
$$R_n= \begin{cases} \emptyset  &\mbox{ if  }n \equiv 0 \bmod 3 \\
\{(n-1,n-1), (n-1,n)\}  &\mbox{ if  }n \equiv 2 \bmod 3\\
\{n-3,n-2,n-1,n\} \times \{n-3,n-2,n-1,n\} &\mbox{ if  }n \equiv 1 \bmod 3.
\end{cases}$$

\begin{lem}\label{sharp} For all $n \ge 3$, if $A$ is a subset of
$[n]
\times [n]$ that contains
$B_n$ then  $\tilde H_{\nu_{2n}}(M(A)) \ne 0$.
\end{lem}

\begin{proof} For $n \equiv 0 \bmod 3$, let $$\rho =  \alpha_{1,1^\prime,2^\prime} \land
\beta_{2,3,3^\prime}
\land
\alpha_{4,4^\prime, 5^\prime}\land \beta_{5,6,6^\prime}\land \dots \land
\alpha_{n-2,(n-2)^\prime, (n-1)^\prime}\land \beta_{n-1,n,n^\prime}, $$
 and for $n \equiv 2 \bmod 3$, let $$\rho =  \alpha_{1,1^\prime,2^\prime} \land
\beta_{2,3,3^\prime}
\land
\alpha_{4,4^\prime, 5^\prime}\land \beta_{5,6,6^\prime}\land \dots \land
\alpha_{n-1,(n-1)^\prime, n^\prime}. $$
In both  cases $\rho$ is a cycle in $C_{\nu_{2n}}(M(A))$, but not a boundary.
Indeed,  if
$\rho
$ were a boundary in  $C_{\nu_{2n}}(M(A))$ then it would be  a boundary in
$C_{\nu_{n,n}}(M_{n,n})$, which would  imply that  all the generators of $\tilde
H_{\nu_{n,n}}(M_{n,n})$ given in Lemmas~\ref{gen0} and \ref{gen1} are boundaries.  This is
impossible since by Theorem~\ref{conj},
 $\tilde H_{\nu_{n,n}}(M_{n,n})
\ne 0$. Hence,  $\tilde H_{\nu_{2n}}(M(A)) \ne 0$.

For $n \equiv 1 \bmod 3$,  let $$\rho =  \alpha_{1,1^\prime,2^\prime} \land
\beta_{2,3,3^\prime}
\land \dots \land
\alpha_{n-6,(n-6)^\prime, (n-5)^\prime}\land \beta_{n-5,n-4,(n-4)^\prime}. $$
By Theorem~\ref{conj} and Lemmas~\ref{gen2} and \ref{gen0}, there is a cycle $\omega$  in
$C_{2}(M(R_n))$ such that   the cycle $\rho \land \omega$ is not a boundary in
$ C_{\nu_{n,n}}(M_{n,n})$.
So  $\rho \land \omega$ is not a boundary in $C_{\nu_{2n}}(M(A))$.
Hence  $\tilde H_{\nu_{2n}}(M(A)) \ne 0$.
\end{proof}

\begin{thm} For  $n \ge 3$,  $\tilde H_{\nu_{2n}}(M(D_n))\ne 0$.
\end{thm}

\begin{proof}
We claim that  an isomorphic copy of $D_n$ contains $B_n$ for all $n \ge 3$
except for $n = 4,7$.  Indeed, if $n \equiv 0,2 \bmod 3$ then the isomorphic copy of $D_n$ is
$$[n] \times [n] \setminus (\{(i,i+2) : i = 1, \dots, n-2\} \cup \{(n-1,1),(n,2)\}). $$
If $n \equiv 1 \bmod 3$ and $n \ge 10$ then the isomorphic copy of $D_n$ is
$$[n] \times [n] \setminus (\{(i,i+4) : i = 1, \dots, n-4\} \cup
\{(i+n-4,i) : i =1,2,3,4 \}).$$  The result now follows from
Lemma~\ref{sharp} and Table 8.1.
\end{proof}

Table 8.1 and the torsion results of Section~\ref{ctorsec} suggest
the following conjecture.
\begin{conj}There exists an integer $n_0 \ge 8$ such that
if $n \ge n_0$ then $\tilde H_{\nu_{2n}}(M(D_n))$ is  an
elementary $3$-group.  Moreover, if $n \ge n_0$ and  $n \equiv 2 \bmod 3 $ then $\tilde
H_{\nu_{2n}}(M(D_n))=
\Z_3$.
\end{conj}

Bj\"orner and Welker's connectivity result is a consequence of a
more general result of Ziegler.  Indeed, Bj\"orner and Welker
\cite{BW} observe that an isomorphic copy of $D_n$ contains the
set $\Gamma(n,2 \nu_{2n}+1 - n)$ described in the following
theorem.

\begin{thm}[Ziegler \cite{Z}] For $0 \le k \le n-1$, let
$$\Gamma(n,k) =\{(i,j) \in [n] \times [n]:  |j-i| \le k \}.$$
  Let $A$ be a subset of $[n] \times [n]$ that contains
$\Gamma(n,2 \nu_{2n} +1- n)$.  Then $M(A)$ is $(\nu_{2n}-1)$-connected.
\end{thm}

Note  that $B_n \subseteq \Gamma(n,2\nu_{2n}+1-n)$ if $n = 6$
 or $n \ge 8$.    It therefore follows from Lemma~\ref{sharp} that  Ziegler's connectivity
bound is sharp for
$n = 6$ and $n \ge 8$.  When $n=3$ or $5$, $M(\Gamma(n,2\nu_{2n}+1 -n))$  is a simplex, which
is contractible.  Hence Ziegler's bound is not sharp in these cases.

\section{Shellability of the $\nu_n$-skeleton of $M_n$} \label{shellsec}

In this section we describe a shelling of the $\nu_n$-skeleton of
$M_n$ along with a discrete Morse function on $M_n$ which is
closely related to our shelling.  We assume that the reader is
familiar with the basic definitions from shellability theory (see
for example \cite{BjWa}) and discrete Morse theory (see
\cite{Fo}). Before presenting our results, we remark that in
\cite{A}, Athanasiadis has shown that the $\nu_n$-skeleton of
$M_n$ is vertex decomposable, which implies that it is shellable.
In light of this fact, we will not provide a proof that our
ordering of the facets of the $\nu_n$-skeleton is in fact a
shelling.

Our shelling and Morse function are determined with use of the
following recursive algorithm, which gives, for any graph $G \in
M_n$, an ordered partition $\rho(G)=(G_1,\ldots,G_r)$ of $G$ into
subgraphs $G_i=(V_i,E_i)$.  We begin with
$G_0=(\emptyset,\emptyset)$. Having defined $G_j$ for all $j<i$,
we define $G_i$ as follows.
\begin{itemize}
\item If $\bigcup_{j<i}V_j=[n]$, stop.
\item If $\bigcup_{j<i}V_j=[n] \setminus \{t\}$, set $G_i=(\{t\},\emptyset)$.
\item If $|\bigcup_{j<i}V_j| < n-1$, let $a,b$ be the two smallest elements of
$[n] \setminus \bigcup_{j<i}V_j$.  Set $V_i=\{a,b\} \cup N_G(a) \cup N_G(b)$ and
define $E_i$ to be the set of all edges in $G$ that have both vertices in $V_i$.
\end{itemize}

For example, if $n=10$ and $E(G)=\{17,38,45\}$ then our algorithm
will give $G_1=(\{1,2,7\},\{17\})$, $G_2=(\{3,4,5,8\},\{38,45\})$,
$G_3=(\{6,9\},\emptyset)$ and $G_4=(\{10\},\emptyset)$.

It follows immediately from the definition of our partition that
$|V_i| \leq 4$ for all $i \in [r]$ and that $|V_i|>1$ if $i<r$.
Moreover, we have $|E_i|=\lfloor \frac{|V_i|}{2} \rfloor$ whenever
$|V_i| \neq 2$.  We now partially order the set of all graphs
$G=(V,E)$ such that $V \subseteq [n]$ by setting $(V,E)\preceq (V^\prime,E^\prime)$ if either
$|V|<|V^\prime|$ or we have
$V=V^\prime=\{i,j\}$ and $E=\emptyset$ while $E^\prime=\{ij\}$.  The
partial order $\preceq$ gives rise to a lexicographic partial
order $\preceq_l$ on $M_n$. That is, if $G,H \in M_n$ with
$\rho(G)=(G_1,\ldots,G_r)$ and $\rho(H)=(H_1,\ldots,H_s)$, we set
$G \preceq_l H$ if either $G_i=H_i$ for all $i \in [r]$ or, for
some $i \leq r$, we have $G_j=H_j$ for all $j<i$ and $G_i \prec
H_i$.

\begin{thm}
Let $F_1<F_2<\ldots <F_t$ be any linear extension of the
restriction of $\preceq_l$ to the set of $\nu_n$-dimensional faces
of $M_n$.  Then $F_1,F_2,\ldots,F_t$ is a shelling of the
$\nu_n$-skeleton of $M_n$. \label{shell}
\end{thm}

To a shelling $F_1,\ldots,F_t$ of any complex $\Delta$, one can
associate a discrete Morse function (actually, many such
functions) as follows.  For each nonhomology facet $F_i$ of the
shelling, let $R_i \subset F_i$ be the restriction face of $F_i$,
that is, the unique minimal new face obtained when $F_i$ is added
to the complex built from $\{F_j:j<i\}$.  The interval $[R_i,F_i]$
in the face poset of $\Delta$ is isomorphic to the face poset of a
simplex (of dimension at least one), and if we fix an isomorphism
between these two posets then any simplicial collapse of the
simplex to a point gives rise to a pairing ${\mathcal M}_i$ of the
faces in $[R_i,F_i]$.  The union of all such pairings ${\mathcal
M}_i$ determines (the gradient flow of) a discrete Morse function
on $\Delta$ whose critical cells are the homology facets of the
given shelling.

A discrete Morse function associated to the shelling of Theorem
\ref{shell} is quite easy to describe. For $G \in M_n$ with
$\rho(G)=(G_1,\ldots,G_r)$, define
\[
\mu(G):=\lp \begin{array}{ll} \infty & \mbox{if no $V_i$ has size
two}, \\ \min\{i:|V_i|=2\} & \mbox{otherwise}. \end{array} \right.
\]
Let $X_n$ be the set of all $G \in M_n$ such that $\mu(G) \neq
\infty$ and $E_{\mu(G)} \neq \emptyset$. For $G \in X_n$, let
$G^-$ be the graph obtained from $G$ by removing the unique edge
in $E_{\mu(G)}$. The next result is straightforward to prove using
standard techniques from discrete Morse theory.

\begin{thm} \label{Morse}
The set $\{(G,G^-):G \in X_n\}$ determines the gradient flow of a
discrete Morse function on $M_n$ whose critical cells are those $G
\in M_n$ such that $\mu(G)=\infty$. \label{morse}
\end{thm}

One can show that the shelling of Theorem \ref{shell} gives rise
to the restriction of the Morse function of Theorem \ref{Morse} to
the $\nu_n$-skeleton of $M_n$.

\section{Bounds on the rank of $\wt{H}_\nu$} \label{ranksec}

In this section we give upper and lower bounds on the rank (that
is, smallest size of a generating set) of $\wt{H}_{\nu_n}(M_n)$
when $n \equiv 0,2 \bmod 3$.  (Note that the case $n \equiv 1 \mod
3$ is settled by Theorem \ref{Btor} and that our lower bound in
the case $n \equiv 0 \bmod 3$ is given in \cite{Bo}.)  We do the
same for $\wt{H}_{\nu_{m,n}}(M_{m,n})$, although we need
conditions on $m,n$ similar to those found in Theorem \ref{ctor}
for the lower bounds.

Set
\[
r_n:=\mbox{rank}(\wt{H}_{\nu_n}(M_n)).
\]
We can get upper bounds on $r_n$ using the Morse function of
Section \ref{shellsec}.  If we let $c_n$ be the size of the set
${\mathcal C}_n$ of graphs $G \in M_n$ with $\nu_n$ edges such
that $\mu(G)=\infty$, then by \cite[Corollary 3.7(i)]{Fo}, we have
\[
r_n \leq c_n.
\]
For $G \in M_n$ with $\rho(G)=(G_1,\ldots,G_r)$, let $\lambda(G)$
be the partition of $n$ such that the number of parts of size $m$
in $\lambda(G)$ is the number of $V_i$ of size $m$.
Straightforward calculation shows that for $G \in M_n$ we have $G
\in {\mathcal C}_n$ if and only if
\[
\lambda(G)=\lp \begin{array}{ll} (3,\ldots,3) & n \equiv 0 \bmod
3, \\ (3,\ldots,3,1) & n \equiv 1 \bmod 3, \\ (4,3,\ldots,3,1) & n
\equiv 2 \bmod 3. \end{array} \right.
\]
Now further calculation gives
\[
c_n=\lp \begin{array}{ll} 2^{n/3} \prod_{j=1}^{n/3}(n-3j+1) & n
\equiv 0 \bmod 3, \\ 2^{(n-1)/3} \prod_{j=1}^{(n-1)/3}(n-3j+1) & n
\equiv 1 \bmod 3, \\ 2^{(n-5)/3}
\sum_{k=1}^{(n-2)/3}\prod_{j=1}^{k}(n-3j+1)\prod_{j=k}^{(n-2)/3}(n-3j)
& n \equiv 2 \bmod 3.
\end{array} \right.
\]
Of course when $n \equiv 1 \bmod 3$ and $n \geq 7$, we know that
$r_n=1$ and our upper bound is both useless and horribly
inaccurate. It turns out that one can improve the upper bound on
$r_n$ in the case $n \equiv 2 \bmod 3$ using the long exact
sequence of Lemma \ref{ontom}. Indeed, if $n \equiv 0 \bmod 3$, the
tail end
\[
\bigoplus_{a,h}\wt{H}_{\nu_n-1}(M_{[n] \setminus
\{1,2,h\}}) \rightarrow \wt{H}_{\nu_n}(M_n)
\rightarrow 0
\]
of the sequence gives
\[
r_n \leq 2(n-2)r_{n-3},
\]
and one simply reobtains the bound $r_n \leq c_n$ using induction.
However, if $n \equiv 2 \bmod 3$ and $n \geq 11$, the tail end of
the sequence is
\[
\bigoplus_{a,h}\wt{H}_{\nu_n-1}(M_{[n] \setminus
\{1,2,h\}}) \rightarrow \wt{H}_{\nu_n}(M_n)
\rightarrow \bigoplus_{i,j}\Z_3 \rightarrow 0,
\]
from which we obtain
\[
r_n \leq 2(n-2)r_{n-3}+(n-2)(n-3).
\]
This recursive formula leads to a somewhat better upper bound than
that given by $c_n$.  However, as we shall see momentarily, all the
bounds we have found so far are so distant from the known lower bounds
on $r_n$ that differences between them are insignificant.  Before going
on to lower bounds, we examine upper bounds for chessboard complexes.
Set
\[
r_{m,n}:= \mbox{rank}(\wt{H}_{\nu_{m,n}}(M_{m,n})).
\]
Using the long exact sequence of Lemma \ref{chessonto} as we used that of Lemma
\ref{ontom} for the matching complexes, we get
\[
r_{m,n} \leq \lp \begin{array}{ll} (m-1)r_{m-2,n-1}+(n-1)r_{m-1,n-2} &
m+n \equiv 0 \bmod 3, \\ (m-1)r_{m-2,n-1}+(n-1)r_{m-1,n-2}+(m-1)(n-1) &
m+n \equiv 2 \bmod 3.
\end{array} \right.
\]
Now we examine lower bounds.  In \cite{Bo}, Bouc gets a lower bound for
$r_n$ when $n \equiv 0 \bmod 3$ using the standard long exact sequence
associated to the pair $(M_n,M_{n-1})$, where we consider $M_{n-1}$ as
the subcomplex of $M_n$ consisting of all matchings in which vertex $n$ is
isolated.  It is straightforward to show that the quotient complex $M_n/M_{n-1}$
has the homotopy type of a wedge of $n-1$ complexes, each homotopy equivalent
to the suspension of $M_{n-2}$, from which it follows that the sequence under
discussion is
\[
\ldots \longrightarrow \wt{H}_t(M_{n-1}) \longrightarrow \wt{H}_t(M_n) \longrightarrow
\bigoplus_{i=1}^{n-1} \wt{H}_{t-1}(M_{n-2}) \longrightarrow \ldots.
\]
When $n \equiv 0 \bmod 3$, the tail end of this sequence is
\[
\wt{H}_{\nu_n}(M_n) \longrightarrow \bigoplus_{i=1}^{n-1}\wt{H}_{\nu_{n-2}}(M_{n-2})
\longrightarrow 0,
\]
from which Bouc obtains
\[
r_n \geq n-1.
\]
When $n \equiv 2 \bmod 3$ and $n \geq 8$, the tail end of the sequence is
\[
\wt{H}_{\nu_n}(M_n) \longrightarrow \bigoplus_{i=1}^{n-1}\wt{H}_{\nu_{n-2}}(M_{n-2})
\longrightarrow \Z_3 \longrightarrow 0,
\]
from which we obtain
\[
r_n \geq (n-1)r_{n-2}-1 \geq (n-1)(n-3)-1.
\]
We can obtain similar results for the chessboard complexes using the long exact sequence
for the pair $(M_{m,n},M_{m-1,n})$, where we consider $M_{m-1,n}$ to be the subcomplex of
$M_{m,n}$ consisting of all matchings in which vertex $m$ is isolated.  We get
\[
r_{m,n} \geq \lp \begin{array}{ll} n & m+n \equiv 0 \bmod 3 \mbox{
and } m \leq n \leq 2m-3,
\\ n(n-1)-1 & m+n \equiv 2 \bmod 3 \mbox{ and } m \leq n \leq 2m-7.
\end{array} \right.
\]
Certainly the distance between the upper and lower bounds we have
provided is unsatisfactory in all cases.

\vspace{.2in}\noindent{\bf Acknowledgements.} We are very grateful to Volkmar Welker,
Frank Heckenbach, and Jean-Guilleaume Dumas for the support they gave us in the use of their
homology software package.  We also thank Xun Dong for pointing out an error in an earlier
version of Section~\ref{shellsec}.

\end{document}